\documentclass[12pt,a4paper,twoside]{article}
\usepackage{amsmath,amsfonts,amssymb,a4,enumerate,epsfig,array,theorem}

\def\qed{\unskip\nobreak\hfil\penalty50\hskip1.75em\null\nobreak\hfil
$\blacksquare$ {\parfillskip=0pt \finalhyphendemerits=0 \par}\goodbreak}

\newfont{\eightrm}{cmr8}
\newfont{\fiverm}{cmr5}

\newcommand\tr{\operatorname{tr}}

\newcommand\clo{\operatorname{cl}}
\newcommand\diag{\operatorname{diag}}
\newcommand\actv{\star}

\newcommand\cA{{\cal A}}
\newcommand\cG{{\cal G}}
\newcommand\cI{{\cal I}}

\newcommand{\VV}{{\cal V}}

\newcommand{\ta}{A^\infty}
\newcommand\nobf{\noindent\bf}
\newcommand\az{\alpha_0}
\newcommand\ai{\alpha_1}
\newcommand\bi{\beta_i}
\newcommand\bo{\beta_o}

\newcommand{\ZZ}{{\mathbb{Z}}}

\newcommand{\RR}{{\mathbb{R}}}
\newcommand{\CC}{{\mathbb{C}}}
\newcommand{\Ss}{{\mathbb{S}}}

\newtheorem{prop}{Proposition}[section]
\newtheorem{coro}[prop]{Corollary}
\newtheorem{theo}[prop]{Theorem}
\newtheorem{defin}[prop]{Definition}
\newtheorem{lemma}[prop]{Lemma}

\begin{document}
\title{Tilings of quadriculated annuli}
\author{Nicolau C. Saldanha and Carlos Tomei}
\maketitle

\begin{abstract}
Tilings of a quadriculated annulus $A$ are counted according to {\em volume}
(in the formal variable $q$) and {\em flux} (in $p$).
We consider algebraic properties of the resulting
generating function $\Phi_A(p,q)$.
For $q = -1$, the non-zero roots in $p$ must be roots of unity
and for $q > 0$, real negative.
\end{abstract}

\section{Introduction}

A {\it quadriculated surface}
\footnote{1991 {\em Mathematics Subject Classification}.
Primary 05B45, Secondary 05A15, 05C50, 05E05.
{\em Keywords and phrases} Quadriculated surfaces, tilings by dominoes,
dimers.}
\footnote{The authors gratefully acknowledge the support of
CNPq, Faperj and Finep (Brazil).
The authors thank the referee for the many suggestions
which have improved the presentation.}
is a juxtaposition of squares along sides forming a surface
(usually with boundary)
such that vertices of squares which lie in the interior
of the surface belong to precisely four squares
and vertices in the boundary belong to finitely many squares.
A simple example of a quadriculated surface is the {\it quadriculated plane},
$\RR^2$ divided into unit squares with vertices in $\ZZ^2$.
Compact quadriculated surfaces have finitely many squares
and will usually be {\it bicolored} and {\it balanced},
i.e., the forming squares are black and white in such a way that
squares with a common side have opposite colors,
and the number of black squares equals the number of white squares.

Let $A$ be a balanced quadriculated annulus:
$A$ is topologically $\Ss^1 \times [0,1]$
and can be embedded in $\RR^2$.
We assume $A$ to be a subset of the plane,
allowing free use of an orientation of $A$
and of its inner and outer boundaries.

Let $T_A$ be the set of domino tilings of $A$.
For $t \in T_A$, we draw inside each domino
a {\it domino arrow} from the center of the black square
to the center of the white square.
Given $t_0, t_1 \in T_A$,
juxtapose the domino arrows from $t_1$ to the reversed domino arrows from $t_0$:
this produces a finite number of oriented simple closed curves.
The {\it flux} $\phi(t_1; t_0)$ of $t_1$ relative to $t_0$
counts such curves with weights $-1$, 0 or 1
given by the winding of each curve in the annulus $A$.
The word {\it flux} has been chosen for its (co)homological connotations;
unfortunately, this use of the term is unrelated to that of
Lieb and Loss (\cite{LL}).
The {\it flux polynomial} of $A$ is
$$\Phi_A^+(p) = \sum_{t \in T_A} p^{\phi(t; t_0)};$$
this is a Laurent polynomial in $p$
which is independent of base tiling $t_0$,
up to multiplication by an integer power of $p$.

\begin{theo}
\label{theo:plus}
All non-zero roots of the flux polynomial $\Phi_A^+(p)$
are real negative numbers.
Furthermore, they are either all equal to $-1$ or all distinct.
\end{theo}
 
When the roots of $\Phi_A^+(p)$ are distinct, $-1$ may or may not be a root.

A tiling $t$ induces a bijection $\beta_t$ from black to white squares.
Again fixing an arbitrary base tiling $t_0 \in T_A$,
let the {\it sign} $\sigma(t; t_0)$ of a tiling $t$
relative to $t_0$ be the sign
of the permutation $\beta_t \circ \beta_{t_0}^{-1}$ on the set of white squares.
Deift and Tomei (\cite{DT}) pointed out that $-1$-counting of tilings
on disks can only yield very few values.

\begin{theo}[Deift-Tomei]
\label{theo:DT}
Let $D$ be a quadriculated disk with a base tiling $t_0$.
Then
\[ \sum_{t \in T_D} \sigma(t; t_0) \]
equals $0$, $1$ or $-1$.
\end{theo}

The {\it signed flux polynomial} of an annulus $A$ is
$$\Phi_A^-(p) = \sum_{t \in T_A} \sigma(t; t_0) p^{\phi(t; t_0)},$$
an expression independent of base tiling
up to sign and multiplication by a power of $p$.
On quadriculated annuli, $-1$-counting is also severely restricted.

\begin{theo}
\label{theo:minus}
All non-zero roots of the signed flux polynomial $\Phi_A^-(p)$
are roots of unity.
\end{theo}

Actually, we study a more general polynomial,
the {\it $q$-flux polynomial} of $A$,
\[ \Phi_A(p,q) = \sum_{t \in T_A} q^{\nu(t;t_0)}p^{\phi(t;t_0)}, \]
which enumerates tilings of $A$ according to {\it volume} $\nu$ and flux.
Counting according to volume (defined in Section 4), i.e.,
$q$-counting of domino tilings,
has been considered previously: in \cite{EKLP}, for example,
Elkies, Kuperberg, Larsen and Propp $q$-count tilings
by dominoes of the {\it Aztec diamond} and
MacMahon (\cite{M}) already $q$-counted lozenge tilings.
The flux polynomials are special cases of the $q$-flux polynomial:
$\Phi_A^+(p) = \Phi_A(p,1)$ and $\Phi_A^-(p) = \Phi_A(p,-1)$.

The main ingredient in the proofs is the fact that $\Phi_A(p,q)$ has
two distinct algebraic interpretations. The first one is rather
familiar: $\Phi_A(p,q)$ is the determinant
of a Kasteleyn-type matrix $M_A(p,q)$ (see Section 4).
The second interpretation requires the construction of another matrix:
for each flux $f$, the coefficient of $p^f$ in $\Phi_A(p,q)$
(a Laurent polynomial in $q$)
is the trace of a {\it connection matrix} obtained as follows (see Section 6).
Tear $A$ along a {\it cut} formed by sides of squares
(as in Figure \ref{fig:tracks})
giving rise to a {\it periodic track segment}:
a (topological) quadriculated rectangle $\Delta$
with a marked pair of opposite sides, the left and right {\it attachments}
(see Figures \ref{fig:segtile} and \ref{fig:nilpot}).
Now, a tiling of $A$ almost induces a tiling on $\Delta$:
the complication is that the
dominoes of the tiling, when drawn so as to fill $\Delta$ completely,
may trespass the attachments.
Our definition of tiling of a track segment
shall allow dominoes to trespass attachments
(again as in Figure \ref{fig:segtile}).
We list, for a given value of the flux $f$,
all possible ways in which dominoes
can trespass the attachment (the {\it indices} of an attachment).
Say there are $N$ such ways: the $N \times N$ connection matrix
$C_{\Delta,f}(q)$ has for $(i,j)$-entry the $q$-counting of tilings of $\Delta$
trespassing left and right attachments
in the way prescribed by indices $i$ and $j$.
Rather unsurprisingly,
the connection matrix depends in an irrelevant fashion on the choice of cut.

Both Kasteleyn and connection matrices are well behaved
with respect to taking the quadriculated $n$-cover $A^n$ of an annulus $A$.
As we shall see in Section 5,
the Kasteleyn matrix
associated to $A^n$ has for determinant a polynomial $\Phi_{A^n}(p,q)$
whose roots, for fixed $q$, are, up to sign, the $n$-th power of the roots of
$\Phi_A(p,q)$ for the same fixed $q$. Thus,
we may study the roots of $\Phi_A(p,q)$  by considering different
$n$-covers of $A$. For connection matrices,
the relationship is even simpler: $C_{\Delta^n,f} = (C_{\Delta,f})^n$
(Section 6).
Here, $\Delta$ is a periodic track segment obtained from $A$ by a cut
and $\Delta^n$ is the {\it juxtaposition} of $n$ copies of $\Delta$,
i.e., the periodic track segment obtained by cutting $A^n$.

The connection matrix $C_{\Delta^n,f}$ has two useful features.
First, the computation of each entry is a counting of tilings
on a disk: indeed, the indices associated to left and right attachments,
obtained by the row and column of an entry,
indicate how to {\it prune} $\Delta$
at the attachments yielding a disk $D$
so that tilings of $\Delta$ with prescribed indices
are tilings of the disk $D$, which are then $q$-counted.
We may then apply Theorem \ref{theo:DT} to infer that,
in the case of $(-1)$-counting,
connection matrices can have very few values for their entries.
By finiteness, there must
be different integers $n_0$ and $n_1$ for which, for all possible values
of the flux $f$, we must have
$(C_{\Delta,f}(-1))^{n_0} = (C_{\Delta,f}(-1))^{n_1}$
and hence $\Phi_{A^{n_0}}(p,-1) = \Phi_{A^{n_1}}(p,-1)$.
Since the roots of $\Phi_{A^n}(p,-1)$ are, up to sign,
$n$-th powers of the roots of $\Phi_A(p,-1)$,
this essentially proves that the non-zero roots of $\Phi_A(p,-1)$
are roots of unity---this is Theorem \ref{theo:minus}.
The second feature of the connection matrix is the simple fact that, by
definition, it has only non-negative entries when $q>0$. This permits
the study of its eigenvalues (and thus, its trace) by making use of
the Perron theorem for non-negative matrices.

In Section 2, we introduce {\it wall-free annuli},
essentially irreducible quadriculated annuli
with respect to the set of tilings. For annuli admitting walls,
the theorems have more stringent statements, with simple proofs.
In this section, we introduce most of the notation used in the paper
and a convenient discrete version of the Gauss-Bonnet theorem.

In Section 3, we review the basic facts about {\it height functions},
algebraic descriptions of tilings
of a simply connected quadriculated region (\cite{T}).
Tilings of an annulus $A$ give rise to {\it periodic} tilings on $\ta$,
which are then studied as height functions in the {\it band} $\ta$.

In Section 4, we define {\it Kasteleyn weights} on $A$.
This is a modification of Kasteleyn's construction (\cite{K}, \cite{K2}).
In more detail, cut $A$ to obtain a periodic track segment $\Delta$
and a {\it Kasteleyn triple} of matrices:
one accounting for interaction between black
and white squares whose common side is not in the cut, and two others
which account for interaction of squares across the cut.
The triple is convenient when studying
the effect of {\it juxtaposition} and {\it closing} of track segments.
The $q$-flux polynomial $\Phi_A(p,q)$ is the
determinant of a matrix obtained in a simple fashion from the Kasteleyn triple.
We build up from three simple homomorphisms---%
{\it flux}, {\it volume} and {\it sign}---%
on the (additive) group of 1-cycles on $A$.
The homomorphisms have simple combinatorial interpretations
and incorporate well the properties
required for obtaining $(p,q)$-counting of tilings out of a determinant.
In Section 5, we show how Kasteleyn triples and the $q$-flux polynomial
behave with respect to $n$-covers of the annulus $A$.

Section 6 describes the connection matrix $C_{\Delta,f}$.
The proof of Theorem \ref{theo:minus} ends the section.
Most of Section 7 is dedicated to the positivity of certain matrix entries
which is needed to use Perron's theorem.
We finally prove the following generalization of Theorem \ref{theo:plus}.

\begin{theo}
\label{theo:qplus}
Let $A$ be a wall-free annulus.
Let $q>0$ be fixed. Then all non-zero roots of the $q$-flux
polynomial $\Phi_A(p,q)$ are distinct, negative numbers.
\end{theo}

There is a natural symbolic dynamical context for the results in this paper
which shall not be explored in detail
(see \cite{KH} for terminology and basic results).
For the band $\ta$, let $T_{\ta,f}$ be the set of tilings of $\ta$
with flux $f$;
$T_{\ta,f}$ is a Cantor set under a natural sequence space topology.
{\it Translation} of tilings by one period corresponds to the shift
in sequence space and is a homeomorphism
from $T_{\ta,f}$ to itself.
With this structure, $T_{\ta,f}$ is a topological Markov chain
with transitive matrix and therefore
topologically mixing and with a dense set of periodic tilings.
Furthermore, translation preserves a certain class of measures
constructed taking into account the combinatorics of tilings.

It seems natural to ask if similar results can be obtained
for tilings by {\it lozenges} of {\it triangulated annuli}
(a lozenge is a juxtaposition of two equilateral triangles
along a common side of length 1).
We indicate two significant difficulties.
First, there is no obvious counterpart to Theorem \ref{theo:DT}.
Besides, the technical requirements to prove eventual compatibility,
indeed, the very concept of wall,
have not, to our knowledge, been sufficiently studied.
Another related problem is the study of other surfaces
such as tori and surfaces of higher genus with boundary.
This must involve polynomials in more variables
to account for the value of the flux across independent cuts.

\section{Geometry and topology of quadriculated \\ annuli}

In the introduction we already presented the concepts of
bicolored and balanced quadriculated surfaces.
We always consider disks and annuli embedded in $\RR^2$
but not necessarily in the quadriculated plane.
Let $T_S$ be the set of domino tilings of a quadriculated surface $S$.
A {\it path} $v_0v_1\ldots v_n$ in $S$ is
a sequence of adjacent vertices.
A path $v_0\ldots v_n$ {\it turns left} (resp. {\it right}) at $v_i$
if $v_{i-1}$, $v_i$ and $v_{i+1}$ are vertices of a square to the left
(resp. right) of both edges $v_{i-1}v_i$ and $v_iv_{i+1}$;
otherwise, if $v_i$ is an interior point in  $S$,
the path {\it follows straight} at $v_i$.
More generally,
the {\it curvature} of the path $\ldots v_{i-1}v_iv_{i+1}\ldots$ at $v_i$ 
is 1, 0 or $-1$ if $v_i$ is in the interior of $S$ and
the path turns left, follows straight or turns right at $v_i$, respectively;
if $v_i$ is a boundary vertex, let $n$ be the number of squares in the angle
$v_{i-1}v_iv_{i+1}$: the curvature at $v_i$ is $\pm(2-n)$,
with the sign depending on the orientation of the angle. 
Clearly, the definition of curvature at interior vertices
can be considered a special case of the more complicated
definition for boundary vertices.
A path is {\it simple} if its vertices are distinct
and {\it closed} if $v_0 = v_n$,
and in this case we write $v_m = v_{m \bmod n}$.

A {\it cut} is a simple path in a quadriculated annulus
joining outer to inner boundary and having exactly one vertex
in each boundary component.
We consider an annulus as a circular juxtaposition
of quadriculated disks with specified {\it attachments}.
We draw cuts and call the disks between consecutive cuts {\it track segments}
since we imagine the annulus to be a circular track for a toy train.
More precisely, a {\it track segment} $\Delta$
is a quadriculated disk whose boundary is divided
counterclockwise into four arcs $\bi$, $\az$, $\bo$ and $\ai$
(for inner boundary, attachment 0, outer boundary and attachment 1).
Attachments will be oriented from $\bo$ to $\bi$:
this coincides with the boundary orientation for $\ai$
but not for $\az$.
A track segment $\Delta$ obtained from a balanced annulus $A$
by a single cut is {\it periodic}.

\begin{figure}[ht]
\begin{center}
\epsfig{height=6cm,file=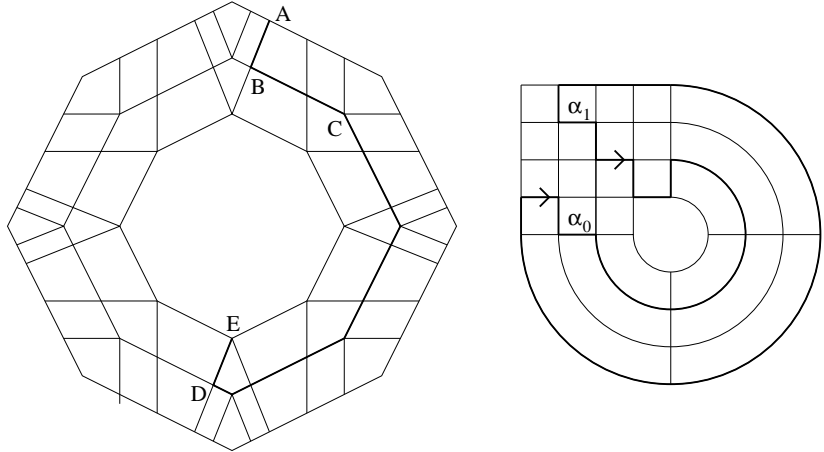}
\end{center}
\label{fig:tracks}
\caption{Closing-up of track segments}
\end{figure}

The {\it shape} of a path $v_0\ldots v_n$ is the sequence
of the $n-1$ curvatures of the path at each $v_i$, $i = 1, \dots, n-1$.
The {\it juxtaposition} of track segments $\Delta$ and $\Delta'$
is natural: attachments $\ai$ and $\az'$ must have the same shape
and their edge-by-edge identification yields the larger track segment
$\Delta\Delta'$ with attachments $\az$ and $\ai'$.
The annulus $A$ is the {\it closing-up} $\clo\Delta$ of the periodic
track segment $\Delta$ and
is obtained by attaching $\ai$ to $\az$;
the resulting cut is said to be {\it induced} by $\Delta$.
Figure \ref{fig:tracks} shows the closing-up of a periodic track segment
and the closing-up of the juxtaposition of two track segments.

The juxtaposition $\Delta^n$ of $n$ copies of a periodic track segment
$\Delta$ under closing-up yields
the {\it $n$-cover} $A^n = \clo(\Delta^n)$ of $A$.
Similarly, the universal cover $\ta$ is also denoted by $\Delta^\ZZ$,
an infinite quadriculated {\it band}.
The canonical projection $\pi: \ta \to A$ and
the {\it translation by one period} $\tau: \ta \to \ta$
are defined as usual.
Notice that $\tau$ takes neighboring vertices to neighboring vertices,
satisfies $\pi \circ \tau = \pi$ and,
for any $p \in \ta$ and any curve $\gamma$ from $p$ to $\tau(p)$,
$\pi \circ \gamma$ goes around $A$ once counterclockwise.
The projection $\pi$ induces a labeling of boundary components
of $\ta$ as inner and outer as well as
orientations on both boundary components;
we continue to call such orientations {\it counterclockwise} or
{\it clockwise}.
The notation $\Delta^\ZZ$ suggests that we see this band as a two-sided
infinite juxtaposition of equal track segments
\[ \ta = \Delta^\ZZ =
\cdots \Delta_{-2} \Delta_{-1} \Delta_0 \Delta_1 \Delta_2 \cdots. \]
The induced cut between $\Delta_i$ and $\Delta_{i+1}$
will be labelled $\xi_{i+\frac{1}{2}}$.
Notice that $\tau$ takes $\Delta_i$ to $\Delta_{i+1}$ and
$\xi_{i+\frac{1}{2}}$ to $\xi_{i+\frac{3}{2}}$.

Let $A$ be a quadriculated annulus, $\xi$ be a cut and $t \in T_\ta$
a tiling of $\ta$.
Draw domino arrows inside each domino in $t$ from black to white squares.
The {\it flux} $\phi(t,\xi)$ of $t$ relative to $\xi$ is
the number of such arrows crossing $\xi$ counterclockwise
minus the number of arrows crossing $\xi$ clockwise.

\begin{lemma}
\label{lemma:fconst}
Let $\Delta$ be a periodic track segment, $A = \clo\Delta$, $t \in T_\ta$
and let $\ta$ have induced cuts $\xi_{n+\frac{1}{2}}$.
The value of $\phi(t; \xi_{n + \frac{1}{2}})$ is the same for all $n \in \ZZ$.
\end{lemma}

{\nobf Proof:}
The numbers of black and white squares in $\Delta_0$ are equal.
Also, the numbers of black and white squares in the union $X$
of all dominoes in $t$ with at least one square in $\Delta_0$
are also equal.
But $x = \phi(t; \xi_{\frac{1}{2}}) - \phi(t, \xi_{-\frac{1}{2}})$ is the number
of white squares in $X - \Delta_0$ minus the number
of black squares in $X - \Delta_0$ and therefore $x = 0$.
The figure shows how each domino contributes to $x$.
\qed

\begin{figure}[ht]
\begin{center}
\epsfig{height=3cm,file=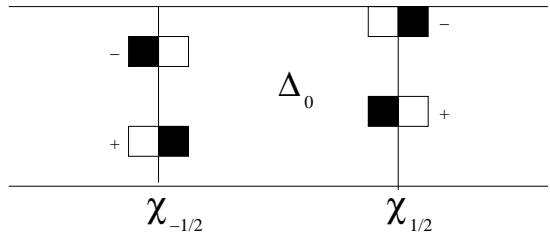}
\end{center}
\caption{Comparing fluxes across cuts}
\label{fig:flux}
\end{figure}

There is a natural concept of {\it translation of a tiling}
$t_{\ta}$ of $\ta$:
$\tau(t_{\ta})$ is formed by the images under $\tau$
of the tiles of $t_{\ta}$.
A tiling $t_{\ta}$ is {\it periodic} if
$\tau(t_{\ta}) = t_{\ta}$.
A tiling $t_A$ of $A$ induces a tiling $\pi^{-1}(t_A)$ of $\ta$:
clearly, a tiling of $\ta$ is induced by a tiling of $A$
if and only if it is periodic.

\begin{figure}[ht]
\begin{center}
\epsfig{height=5cm,file=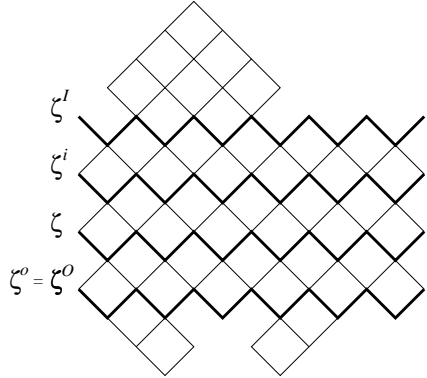} 
\end{center}
\caption{Walls are zig-zags}
\label{fig:zigzag}
\end{figure}

A path is a {\it zig-zag} if it turns alternately left and right
(Figure \ref{fig:zigzag} shows four zig-zags).
Certain quadriculated annuli $A$ admit a simple closed zig-zag $\zeta$.
Keeping in mind that $A$ is embedded in the plane,
$\zeta$ splits the plane in an inner and an outer region
and splits in turn $A$ as a disjoint union of connected components.
Usually these components consist of two annuli,
as will be further discussed in Proposition \ref{prop:closedzigzag}
If, however, the zig-zag hits the boundary of $A$,
the number of components may increase but they will always be disks and annuli;
the zig-zag may also coincide with one of the two boundary components of $A$.
A simple closed zig-zag $\zeta$ in a balanced annulus $A$ is a {\sl wall}
if the numbers of black and white squares 
of $A$ in the bounded component of $\RR^2 - \zeta$ are equal.

\begin{prop}
Let $A$ be a balanced quadriculated annulus with a wall $\zeta$.
No domino in any tiling of $A$ may trespass $\zeta$.
\end{prop}

{\nobf Proof: }
Let $X$ be the set of squares of $A$ in the bounded component
of $\RR^2 - \zeta$ and $Y \supseteq X$ be the union of all dominoes
with at least one square in $X$.
Both $X$ and $Y$ are balanced but all squares in $Y - X$ are of the same color.
Thus $Y - X$ is empty.
\qed

An annulus $A$ with no walls is {\it wall-free}.
From the proposition above, if a wall $\zeta$
splits $A$ as a union of two annuli $A_1$ and $A_2$ then
the flux polynomials satisfy
$\Phi^+_A = \Phi^+_{A_1} \Phi^+_{A_2}$ and
$\Phi^-_A = \Phi^-_{A_1} \Phi^-_{A_2}$.

If a wall $\zeta$ does not touch the inner (resp. outer)
boundary of $A$, then $\zeta$ may be {\it shifted inwards}
(resp. {\it outwards}) to obtain another parallel wall $\zeta^i$
(resp. $\zeta^o$);
we may of course repeat the procedure as often as necessary
to obtain parallel walls $\zeta^I$ and $\zeta^O$
which touch the inner and outer boundary,
respectively (see Figure \ref{fig:zigzag}).
This decomposes $A$ as a union of {\it ladders}
(annuli between consecutive parallel walls) and disks
such that any tiling of $A$ restricts to a tiling of each ladder and disk.

\begin{figure}[ht]
\begin{center}
\epsfig{height=20mm,file=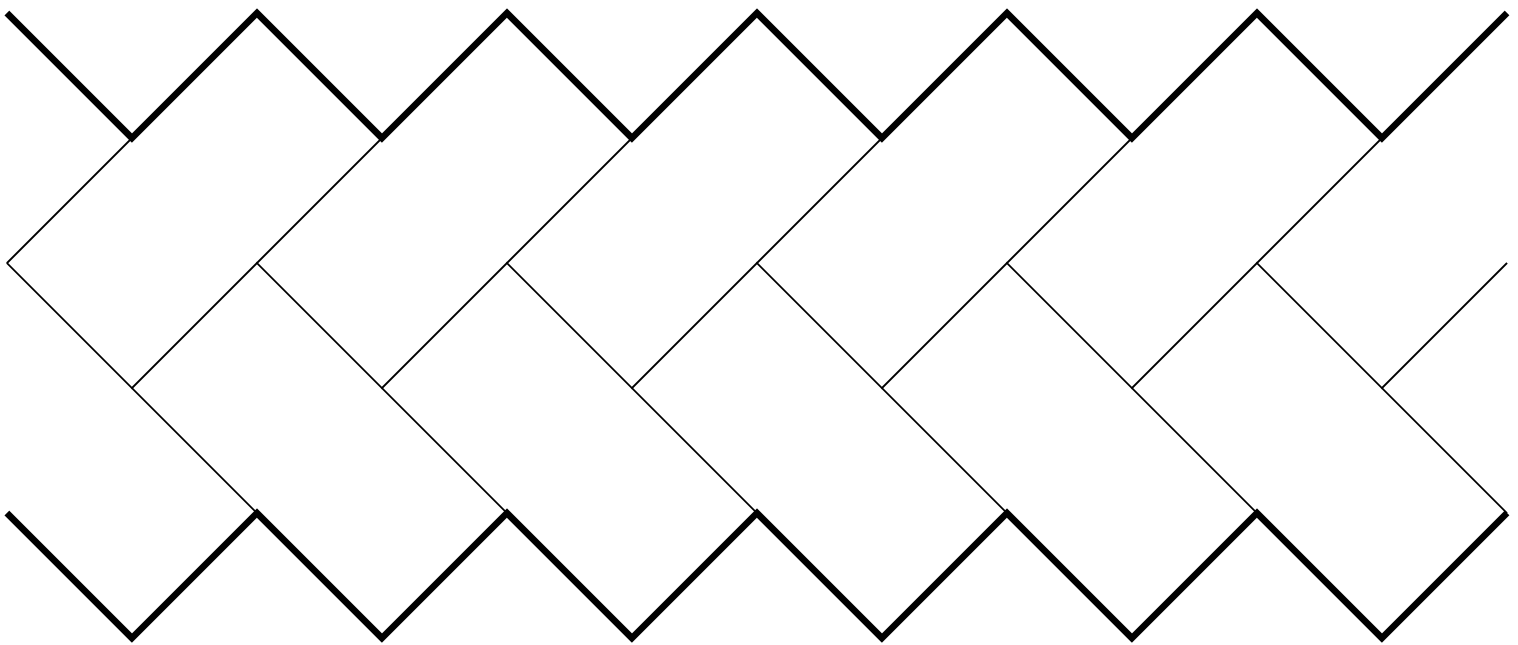}
\quad \epsfig{height=20mm,file=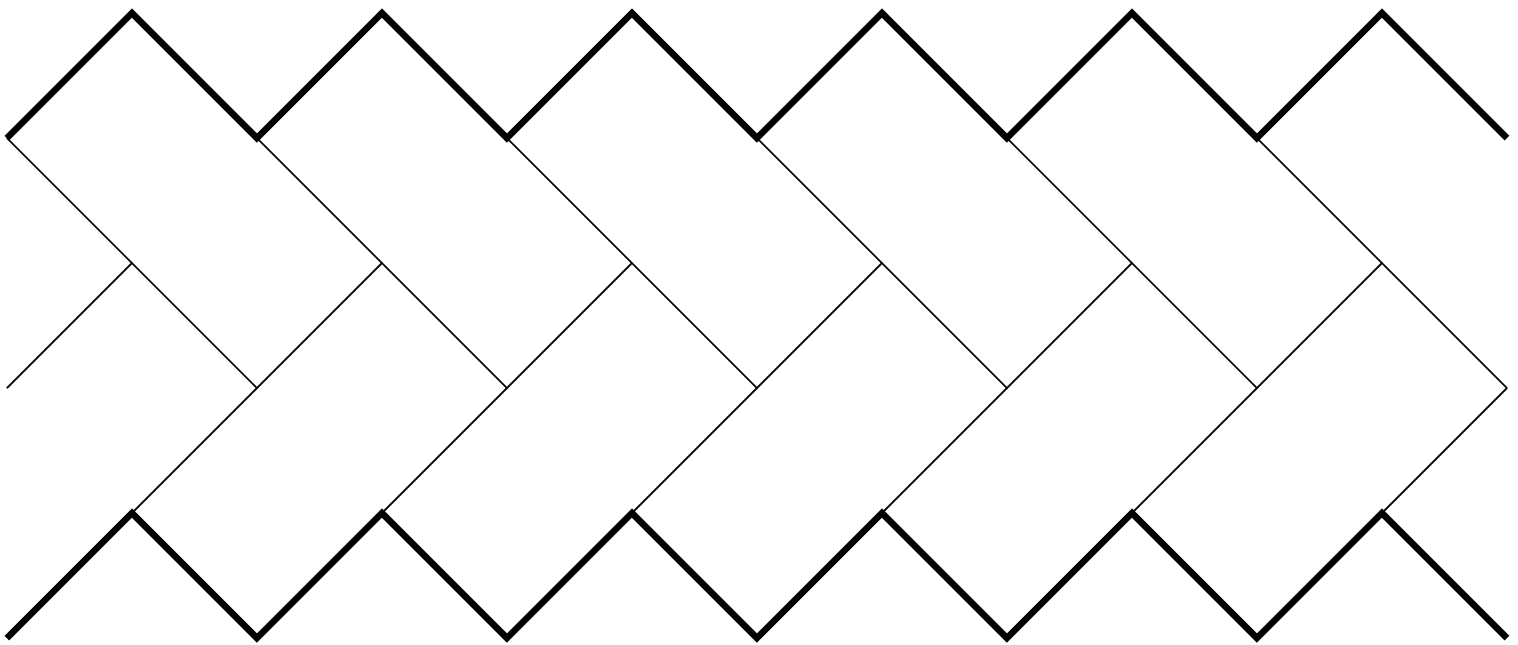}
\end{center}
\begin{center}
\epsfig{height=20mm,file=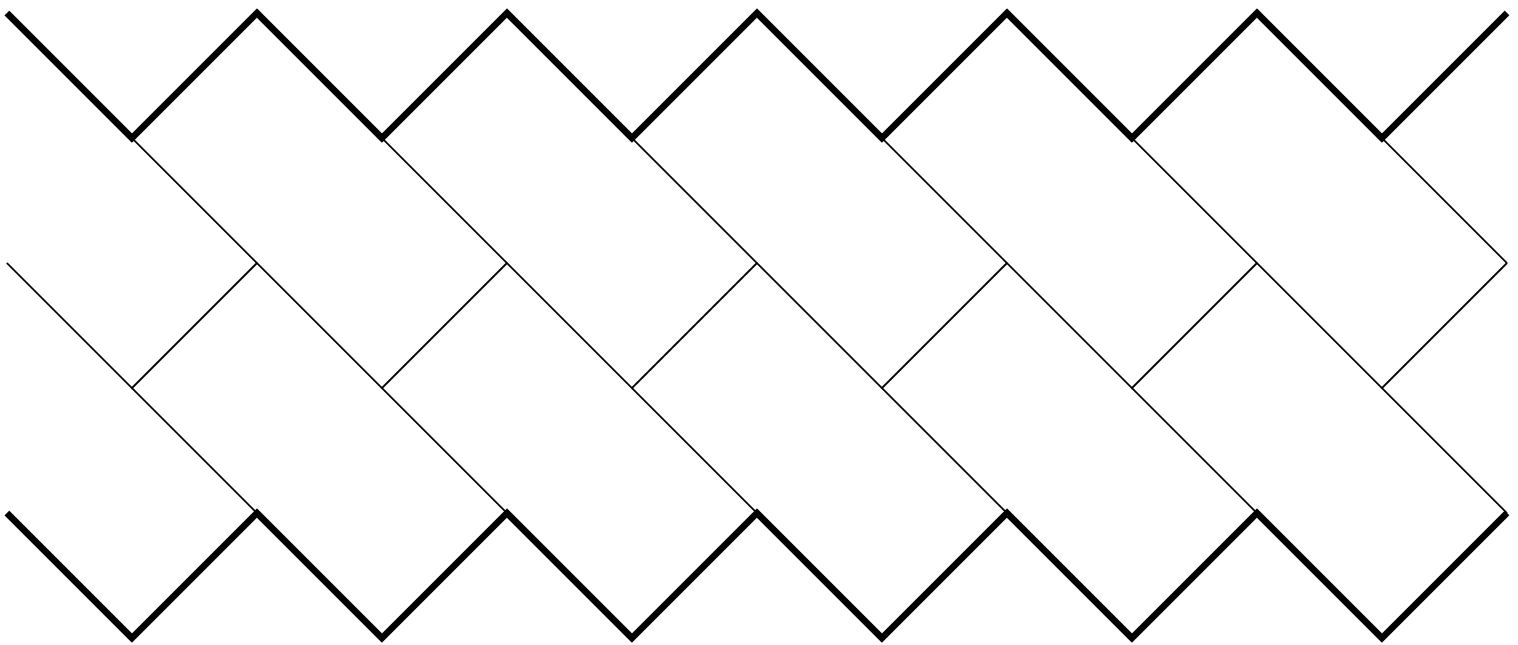}
\quad \epsfig{height=20mm,file=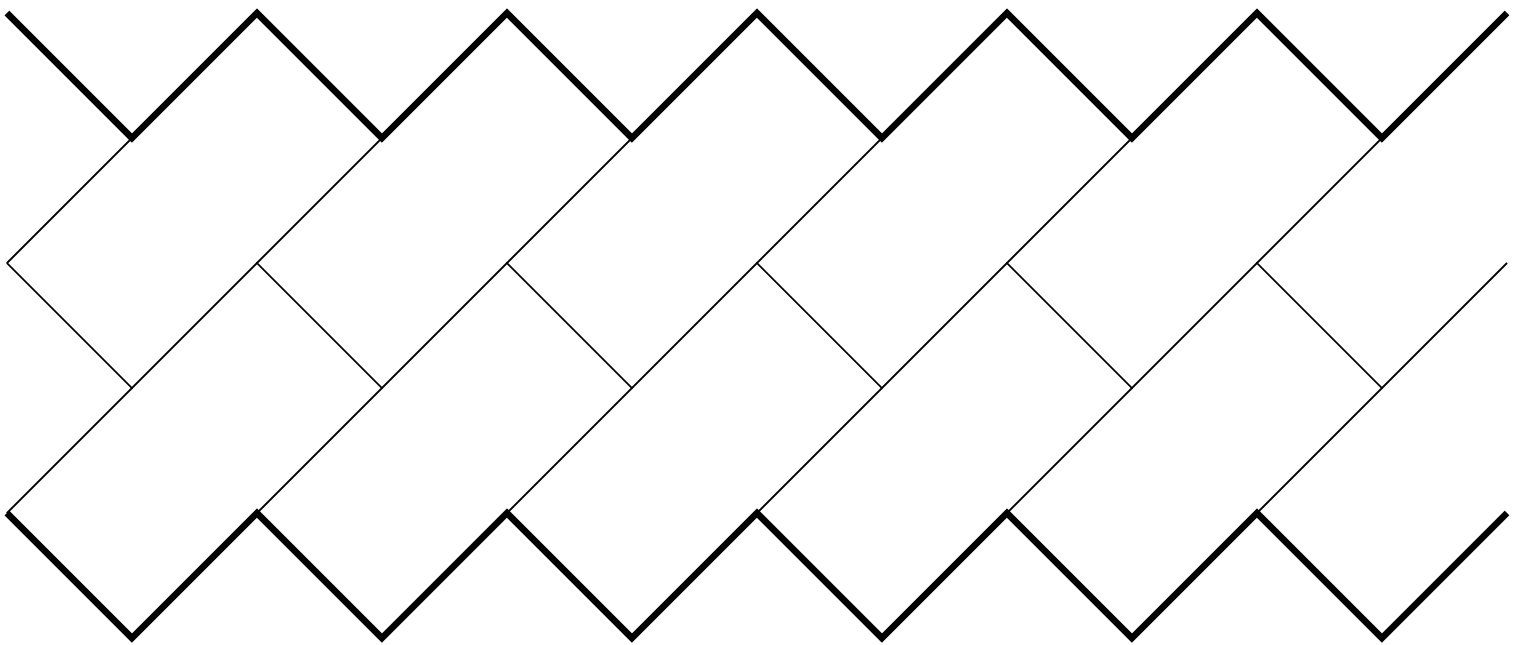}
\end{center}
\caption{Ladders}
\label{fig:ladder}
\end{figure}

A {\it ladder in a tiling} is a sequence of parallel dominoes side by side
such that two neighboring dominoes always touch
along an edge of the longer side,
each domino in the ladder has two neighbors in it
and these two neighbors touch the domino at different squares.
Ladders are then either periodic (in an annulus) or bi-infinite (in a band).
Periodic ladders in a tiling are always contained between two walls.
The tilings in Figure \ref{fig:ladder}
(which may be interpreted as either a band or a cut annulus)
are the four tilings of the surface consisting of two ladders.
Notice that inner and outer boundaries of the annulus
in Figure \ref{fig:ladder} are walls.

\begin{theo}
\label{theo:walls}
Let $A$ be a quadriculated annulus admitting $n+1$ parallel walls, $n \ge 0$.
Then
\[ \Phi_A^+(p) = C^+ (p+1)^n,\qquad \Phi_A^-(p) = C^- (p-1)^n, \]
for constants $C^+$ and $C^-$.
In particular, all roots of $\Phi_A^+(p)$ (resp. $\Phi_A^-(p)$)
are equal to $-1$ (resp. $1$).
\end{theo}

{\nobf Proof:}
Since tilings of $A$ do not trespass walls,
$p$-counting tilings of $A$ amounts to $p$-counting tilings of
connected components of the complement of walls in $A$
(i.e., of disks and ladders) and multiplying the results.
If $L$ is a ladder, $\Phi^+_{L} = p+1$ and $\Phi^-_{L} = p-1$,
up to sign and multiplication by an integer power of $p$.
The restriction of a tiling of $A$ to a connected component which is a disk
bears no relationship with the flux of the tiling:
disks therefore only contribute with multiplicative constants.
\qed

This proves Theorems \ref{theo:plus} and \ref{theo:minus}
and gives a fairly complete description of the set of tilings
for annuli which admit walls.
From now on we may consider wall-free annuli but this hypothesis
will only be crucial in the final section.

The total curvature of a closed path is the sum of the curvatures
at its vertices.
Interior vertices on a zig-zag alternate between curvatures 1 and $-1$
and closed zig-zags have total curvature 0.
Paths with curvature 0 at all interior vertices
are a natural quadriculated version of geodesics
and many of the results for zig-zags hold true for geodesics.
In this paper we have no use for this kind of path
and we drop the subject altogether.

Let $v \in \partial S$ be a boundary vertex of
an oriented quadriculated surface $S$:
the {\it curvature of the boundary} at $v$ is
$2 - n_v$ if there are $n_v$ squares adjacent to $v$.
This coincides with the previous definition of curvature
of the boundary component once it is given the orientation
which leaves the surface to the left.
The {\it total boundary curvature} of a compact quadriculated surface
is the sum of the curvatures of each boundary vertex.
The next proposition is the quadriculated version
of the Gauss-Bonnet theorem;
notice that the usual $2\pi$ is replaced by 4.

\begin{prop}
The total curvature of the boundary of an oriented
compact qua\-dri\-cu\-la\-ted surface $S$ is $4 \chi(S)$.
\end{prop}

The Euler characteristic of $S$ is $\chi(S) = F - E + V$,
the number of faces (squares) minus the number of
edges (sides of squares) plus the number of
(interior and boundary) vertices.

{\nobf Proof }
Let $V_i$ be the number of interior vertices.
Counting vertices face by face and again edge by edge yields
\[ V = V_i + \sum_{v \in \partial S} 1, \quad
4F = 4 V_i + \sum_{v \in \partial S} n_v, \quad
2E = 4 V_i + \sum_{v \in \partial S} (n_v + 1).\]
These three identities imply
\[ 4(F - E + V) = \sum_{v \in \partial S} (2 - n_v), \]
the desired result.
\qed

Thus, since annuli have Euler characteristic 0,
the curvature of the outer and inner boundary components add to 0.
The curvature of the outer boundary component may assume any integer value.

\begin{prop}
\label{prop:closedzigzag}
Let $A$ be a quadriculated annulus embedded in $\RR^2$.
\begin{enumerate}[\rm (a)]
\item{If a closed zig-zag $\zeta$ exists in $A$ then
all zig-zags in $A$ are simple,
each boundary component of $A$ has curvature $0$ and
another maximal zig-zag intersecting $\zeta$
joins both boundary components of $A$.}
\item{If a self-intersecting maximal zig-zag $\zeta$ exists
then the curvatures of the boundary components are $1$ and $-1$
and the endpoints of $\zeta$ are both on the boundary component
with positive curvature.}
\item{For arbitrary $A$, there exists a simple maximal zig-zag 
joining both boundary components.}
\end{enumerate}
\end{prop}

{\nobf Proof: }
Closed zig-zags have total curvature zero.
Also, zig-zags intersect by sharing an edge.
For a self-intersecting zig-zag $\zeta$,
construct a simple closed path $\hat\zeta$ of minimal length
by starting at an edge used twice
and following the zig-zag until returning to that edge.
The path $\hat\zeta$ has total curvature $\pm 1$
and from Gauss-Bonnet cannot bound a subdisk of $A$.
For the same reason, simple closed zig-zags cannot bound subdisks of $A$.
Thus, closed zig-zags are simple and wind once around $A$.

Let $\beta_i$ and $\beta_o$ be the inner and outer boundary components of $A$.
If there exists a closed zig-zag $\zeta$ then,
applying the technique of proof of 
Gauss-Bonnet for the region bounded by $\zeta$ and $\beta_o$,
we see that the curvatures of $\zeta$ and $\beta_o$ add to zero.
Since the curvature of $\zeta$ is zero,
the curvature of both $\beta_o$ and $\beta_i$ are then also zero.
Similarly, if a self-intersecting zig-zag exists then
the curvatures of $\beta_o$ and $\beta_i$ are $\pm 1$.

To finish item (a), take a closed zig-zag $\zeta$.
From an edge of $\zeta$ draw the other possible zig-zag $\zeta^\perp$:
from Gauss-Bonnet, $\zeta^\perp$ cannot intersect $\zeta$ a second time
and it must eventually hit both boundary components.

For item (b), use Gauss-Bonnet now to show that
the self-intersecting zig-zag $\zeta$ cannot self-intersect twice
and its endpoints must both be on the boundary component
of positive curvature.

For item (c), we make use of {\it wiggling} paths:
a path $v_0\ldots v_n$ is {\it wiggling} if its curvature is $\pm 1$
at every vertex $v_i$, $0 < i < n$.
In other words, a path is wiggling if it always has squares of the same color
to its left (or right).
An arbitrary path $v_0\ldots v_n$ may be transformed in a wiggling path
having, say, white to its left and with the same endpoints
by substituting each edge $v_iv_{i+1}$ of the path with black to its left
by $v_iv_{i'}v_{i''}v_{i+1}$, where these are the four vertices of a square.
An edge $v_iv_{i+1}$, $0< i < n-1$ in a wiggling path is {\it zig-zig}
if the curvatures at $v_i$ and $v_{i+1}$ are equal.
Clearly, zig-zags are wiggling paths with no zig-zig edges.

\begin{figure}[ht]
\begin{center}
\epsfig{height=3cm,file=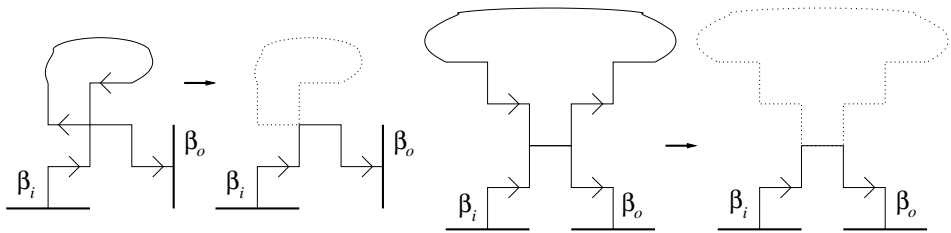}
\end{center}
\caption{Wiggling paths may be taken simple}
\label{fig:wiggling}
\end{figure}

Consider now the class of wiggling paths of minimal length joining
both boundary components.
Such paths are simple, as Figure \ref{fig:wiggling} should convince the reader.
If there is a zig-zag in this class, we found the required zig-zag cut;
this is not always the case so in some cases we have to work more.

\begin{figure}[ht]
\begin{center}
\epsfig{height=3cm,file=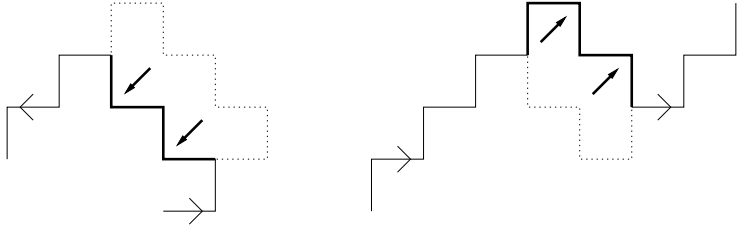}
\end{center}
\caption{Improving wiggling paths}
\label{fig:zigzig}
\end{figure}

What is true instead is that if the class contains no zig-zag
it contains a path with a single zig-zig edge $e$.
Indeed, choose in the class a path such that the length
from the inner boundary to the first zig-zig is maximal.
If the path contained a second zig-zig of identical curvature,
we could construct a shorter wiggling path, a contradiction
(see Figure \ref{fig:zigzig}).
If, on the other hand, the path contained a second zig-zig
of opposite curvature, we could exhibit another wiggling path which
is either shorter or has the same length but a longer initial segment,
again a contradiction (see Figure \ref{fig:zigzig}).
Both constructions indicated in the figure are well defined
in the sense that the existence of the required extra vertices and edges 
follows from the fact that the stretch of path being modified is
initially interior.

We thus have two simple zig-zag segments $\zeta_i$ and $\zeta_o$,
starting respectively from $\beta_i$ and $\beta_o$
with a common edge $e$.
Continue $\zeta_o$ as a zig-zag beyond $e$ until it hits the boundary
or self-intersects.
If it hits the inner boundary we are done.
If it hits the outer boundary,
then $\zeta_o$ together with a piece of $\beta_o$
bound a quadriculated disk into which the zig-zag continuation of $\zeta_i$
enters  at the edge $e$.
Thus, from Gauss-Bonnet, the zig-zag continuation of $\zeta_i$
eventually hits $\beta_o$
without self-intersecting or intersecting the extended $\zeta_o$:
$\zeta_i$ is then the required zig-zag.
Finally, if $\zeta_o$ self-intersects, we are in the situation of item (b)
and $\beta_o$ has curvature $+1$.
The zig-zag continuation of $\zeta_i$ cannot self-intersect (from (b))
and cannot intersect $\zeta_o$ a second time (from Gauss-Bonnet)
and thus must eventually hit $\beta_o$ and we are done.
\qed

\begin{coro}
Let $A$ be a balanced quadriculated annulus and $n$ a positive integer.
Then $A$ admits closed zig-zags if and only if $A^n$ does.
Furthermore, $A$ is wall-free if and only if $A^n$ is wall-free.
\end{coro}

{\nobf Proof:}
A closed zig-zag $\zeta$ in $A$ clearly lifts to another in $A^n$.
Conversely, the projection of a closed zig-zag $\zeta_n$ in $A^n$ to $A$
can not self intersect (by item (a) of the previous proposition)
and is therefore a closed zig-zag $\zeta$ in $A$.
Also, let $d$ be the number of white squares minus the number of black squares
in the bounded connected component of $\RR^2 - \zeta$:
the corresponding difference for $\zeta_n$ is $nd$.
Thus, $\zeta$ is a wall if and only if $\zeta_n$ is.
\qed

\section{Height functions}

We recall the basic properties of {\it height functions}
for a connected and simply connected quadriculated surface $S$
with boundary (\cite{T}; see also \cite{STCR} for the minor adjustments
required for surfaces with no boundary).
Let $\VV_S$ be the set of vertices of squares in $S$
and $\VV_{\partial S} \ne \emptyset$ be the subset of vertices
on the boundary of $S$.
Call an arbitrary but fixed vertex $v_b \in \VV_{\partial S}$
the {\it base vertex}.
For a tiling $t \in T_S$,
construct the height function $\theta: \VV_S \to \ZZ$ as follows:
\begin{enumerate}[(a)]
\item{$\theta(v_b) = 0$;}
\item{if there is a white (resp. black) square to the left
of the oriented edge $v_0v_1$ not covered by a domino of $t$ then
$\theta(v_1) - \theta(v_0) = 1$ (resp. $-1$).}
\end{enumerate}
Clearly, these rules are consistent when constructing $\theta$
along the path surrounding a domino.
This local consistency, together with the fact that $S$
is connected and simply connected,
yield uniqueness and global consistency of $\theta$.
Notice that all height functions coincide on the boundary
and that changing the base vertex alters $\theta$
by an additive integer constant.
Height functions admit the intrinsic characterization below,
given in strictly local terms
(for a proof and some applications, see \cite{STCR}).

\begin{prop}
\label{prop:height}
Let $S$ be a connected and simply connected quadriculated surface
with a base vertex $v_b \in \VV_{\partial S}$.
Let $\theta$ be an integer valued function defined on $\VV_S$.
Then $\theta$ is the height function of some tiling $t$
(relative to the base point $v_b$)
if and only if the following conditions hold:
\begin{enumerate}[\rm (a)]
\item{$\theta(v_b) = 0$;}
\item{if the oriented edge $v_0v_1$
is in the boundary of $S$ and there is a white (resp. black)
square to its left then $\theta(v_1) - \theta(v_0) = 1$ (resp. $-1$);}
\item{if the oriented edge $v_0v_1$
is in the interior of $S$ and there is a white (resp. black)
square to its left then $\theta(v_1) - \theta(v_0) = 1$ or $-3$
(resp. $-1$ or $3$).}
\end{enumerate}
\end{prop}

Recall that all height functions $\theta$ with base point $v_b$
equal 0 at $v_b$.
From the local characterization of height functions,
given a point $v \in \VV_S$ at a distance $\ell$
(measured along edges) from $v_b$,
$|\theta(v)| \le 3 \ell$.
Thus, even for infinite connected and simply connected quadriculated surfaces,
height functions are locally bounded in the sense
that given any finite subset $X$ of $\VV_S$
the restrictions of all height functions to $X$ 
are bounded by a constant $C_X$.
In particular, the maximum or minimum of any (finite or infinite)
set of height functions is well defined and,
again from the local characterization, is a height function.
The set of height functions, and therefore the set of tilings,
thus form a lattice.

Since the quadriculated annulus $A$ is not simply connected,
the construction of a height function as above is not necessarily
globally consistent.
In \cite{STCR}, this difficulty is addressed by the construction
of {\it height sections}.
Here, we follow the simpler alternative of considering
tilings in the band $\ta$, which correspond
to height functions from $\VV_{\ta}$ to $\ZZ$.
For convenience, we assume the base point $v_b$
to belong to the outer boundary of $\ta$.

\begin{lemma}
\label{lemma:heightflux}
For a fixed reference vertex $v_r$ on the inner boundary of $\ta$,
$\theta(v_r) = 4 \phi(t; \xi_{\frac{1}{2}}) + c$
for all tilings $t$ of $\ta$ and some integer constant $c$.
\end{lemma}

{\nobf Proof:}
If the base vertex $v_b$ and the reference vertex $v_r$
are the two ends of $\xi_{\frac{1}{2}}$ then $\theta(v_r)$
can be computed from $\theta(v_b) = 0$
using the definition of height function by following $\xi_{\frac{1}{2}}$.
More precisely, if $v_{i}$ and $v_{i+1}$ are consecutive vertices
on $\xi_{\frac{1}{2}}$ then $\theta(v_{i+1}) - \theta(v_i)$ equals
\begin{itemize}
\item{$3$ (resp. $-3$) if a domino crosses $\xi_{\frac{1}{2}}$ counterclockwise
(resp. clockwise) at the edge $v_iv_{i+1}$;}
\item{$1$ (resp. $-1$) if no domino crosses $\xi_{\frac{1}{2}}$ at $v_iv_{i+1}$
and the square to the left of $v_iv_{i+1}$ is white (resp. black).}
\end{itemize}
Let $n_{3}$, $n_{-3}$, $n_1$ and $n_{-1}$ be the number of edges
in $\xi_{\frac{1}{2}}$ in each situation above.
Clearly,  \[ \theta(v_r) = -3  n_{-3} - n_{-1} + n_1 + 3 n_3 =
(n_1 + n_{-3} - n_{-1} - n_{3}) + 4 (n_3 - n_{-3}). \]
By definition, $\phi(t; \xi_{\frac{1}{2}}) = n_3 - n_{-3}$.
But $n_1 + n_{-3}$ (resp. $n_{-1} + n_3$) is the number of edges
in the cut with a white (resp. black) square to the left
and therefore $c =  n_1 + n_{-3} - n_{-1} - n_{3}$
does not depend on the tiling $t$.

In the general case when $v_b$ or $v_r$ are not ends of $\xi_{\frac{1}{2}}$,
one has to account for the variations of $\theta$ along stretches of boundary
from ends to $v_b$ or $v_r$ which again are independent of tiling.
\qed

Notice that a height function $\theta$ corresponding to
a periodic tiling of $\ta$ does not usually satisfy
$(\theta \circ \tau)(v) = \theta(v)$ but instead
$(\theta \circ \tau)(v) = \theta(v) + c'$ for some integer constant $c'$;
$c'$ depends on $A$ only and not on the tiling $t$.
With a slight abuse, we still call such height functions {\it periodic}.
Call $\theta_{f,\max}$ (resp. $\theta_{f,\min}$)
the maximum (resp. minimum) of all height functions with a given flux $f$.
These extremal height functions are periodic.
Indeed, from the lattice structure on the set of height functions
and maximality of $\theta_{f,\max}$, it is easy to see that
\[ \theta_{f,\max} = \max(\theta_{f,\max},\theta_{f,\max} \circ \tau - c')
\theta_{f,\max} \circ \tau - c', \]
where $c' = \theta_{f,\max}(\tau(v_b))$.

The following result is an application of height functions
to tilings of annuli.

\begin{theo}
\label{theo:maxflocut}
Let $A$ be a quadriculated annulus 
and let $f_{\max}$ and $f_{\min}$ be the maximum and minimum values
of the flux among all tilings of $A$
relative to an arbitrary cut $\xi$.
Then there exist cuts $\xi_{\max}$ and $\xi_{\min}$
such that no tiling of flux $f_{\max}$ trespasses $\xi_{\max}$
and no tiling of flux $f_{\min}$ trespasses $\xi_{\min}$.
\end{theo}

Like walls (in Theorem \ref{theo:walls}),
the non-trespassable cuts above allow for a multiplicative principle
to be applied to the counting of tilings of extremal flux.
Indeed, let $a^{[n]}_k$ be the number of tilings of the covering space $A^n$
with flux $k + f_{\min}$, $k = 0, \ldots, f_{\max} - f_{\min}$;
we have
$a^{[n]}_{f_{\max} - f_{\min}} = (a^{[1]}_{f_{\max} - f_{\min}})^n$ and
$a^{[n]}_0 = (a^{[1]}_0)^n$.
More general relations among the $a^{[n]}_k$ will be presented
in Proposition \ref{prop:ank}.

{\nobf Proof:}
We only prove the existence of $\xi_{\max}$.
We are searching for a cut on which
all height functions with flux $f_{\max}$ coincide:
this happens if and only if
$\theta_{f_{\max},\min}(v) = \theta_{f_{\max},\max}(v)$
for $v$ on the cut.
Let
\[ X = \left\{ { v \in \VV_{\ta} | 
\theta_{f_{\max},\min}(v) = \theta_{f_{\max},\max}(v) }  \right\} \]
be the set of vertices where all height functions
of flux $f_{\max}$ coincide.
This set contains the outer boundary since $v_b$ belongs to it
and, from Lemma \ref{lemma:heightflux}, also contains the inner boundary.
Let $Y$ be the connected component of $X$
containing the outer boundary of $\ta$
(i.e., the maximal subset of $X$ containing the outer boundary
and such that if two vertices $v_1, v_2 \in X$  are adjacent
and $v_1 \in Y$ then $v_2 \in Y$):
we must prove that $Y$ also contains the inner boundary. Set
\[ \theta(v) =
\begin{cases}
\theta_{f_{\max},\min}(v),& v \in Y, \\
\theta_{f_{\max},\min}(v) + 4,& v \in \ta - Y.
\end{cases}
\]
The function $\theta$ is a height function:
we check the local conditions in Proposition \ref{prop:height}.
The function $\theta$ clearly satisfies items (a) and (b).
As to item (c), if $v_1$ and $v_2$ are adjacent and
either both belong to $Y$ or neither belongs to $Y$
then $\theta(v_1) - \theta(v_2) =
\theta_{f_{\max},\min}(v_1) - \theta_{f_{\max},\min}(v_2)$.
If $v_1 \in Y$ and $v_2 \not\in Y$ then $v_2 \not\in X$ and we have,
from the adjacency of $v_1$ and $v_2$,
$\theta_{f_{\max},\max}(v_2) - \theta_{f_{\max},\min}(v_2) = 4$
and therefore
$\theta(v_1) - \theta(v_2) =
\theta_{f_{\max},\max}(v_1) - \theta_{f_{\max},\max}(v_2)$.
Hence in all cases $\theta(v_1) - \theta(v_2)$ coincides
with the difference in values at $v_1$ and $v_2$ of either
$\theta_{f_{\max},\min}$ or $\theta_{f_{\max},\max}$,
genuine height functions, thus showing that $\theta$ satisfies (c).
If $Y$ contains the inner boundary then $\phi(\theta,\xi) = f_{\max}$
but otherwise $\phi(\theta,\xi) = f_{\max} + 1$,
contradicting the maximality of $f_{\max}$.
\qed

Notice that $f_{\max} - f_{\min}$ is the number of non-zero roots of
$\Phi_A^+(p)$.
Also, from the definition of flow relative to a cut,
$f_{\max} - f_{\min}$ is less than or equal to the length
of the shortest cut in $A$.

\section{Matrices and flux polynomials}

Kasteleyn (\cite{K}) showed how to count tilings 
of quadriculated disks by computing a determinant;
this idea has been extended by many authors,
in particular by producing a matrix whose determinant $q$-counts
tilings of quadriculated disks (\cite{EKLP}).
A variation of these ideas, to be detailed in Proposition \ref{prop:qflux},
obtains the {\it $q$-flux polynomial} $\Phi_A(p,q)$ 
for a quadriculated annulus $A$.
The $q$-flux polynomial will extend the flux polynomials
defined in the introduction:
$\Phi_A^+(p) = \Phi_A(p,1)$ and $\Phi_A^-(p) = \Phi_A(p,-1)$.

Let $\Delta$ be a periodic track segment (as defined in Section 2),
$A = \clo\Delta$ and $\xi$ the induced cut in $A$.
The {\it adjacency graph} $G_A$ of the quadriculated annulus $A$
is the oriented graph whose vertices are squares of $A$
and whose edges connect squares sharing a common side,
always oriented from black to white.
The domino arrows of a tiling are edges in $G_A$.
An edge {\it crosses $\xi$ counterclockwise} (resp. {\it clockwise})
if it leaves through the attachment $\az$ (resp. $\ai$)
and enters through $\ai$ (resp. $\az$).
Let $C_1(A)$ be the set of formal sums (with integer coefficients)
of edges of $G_A$.
Also, let $Z_1(A) \subseteq C_1(A)$ be the additive subgroup
of formal sums such that at each vertex of $G_A$
the sum of the coefficients of the adjacent edges is zero.
For example, the difference $c = t_1 - t_2$ of two tilings of
a quadriculated surface $A$ is an element of $Z_1(A)$.
Following the homological terminology,
the elements of $Z_1(A)$ are {\it $1$-cycles}.

\begin{figure}[ht]
\begin{center}
\epsfig{height=6cm,file=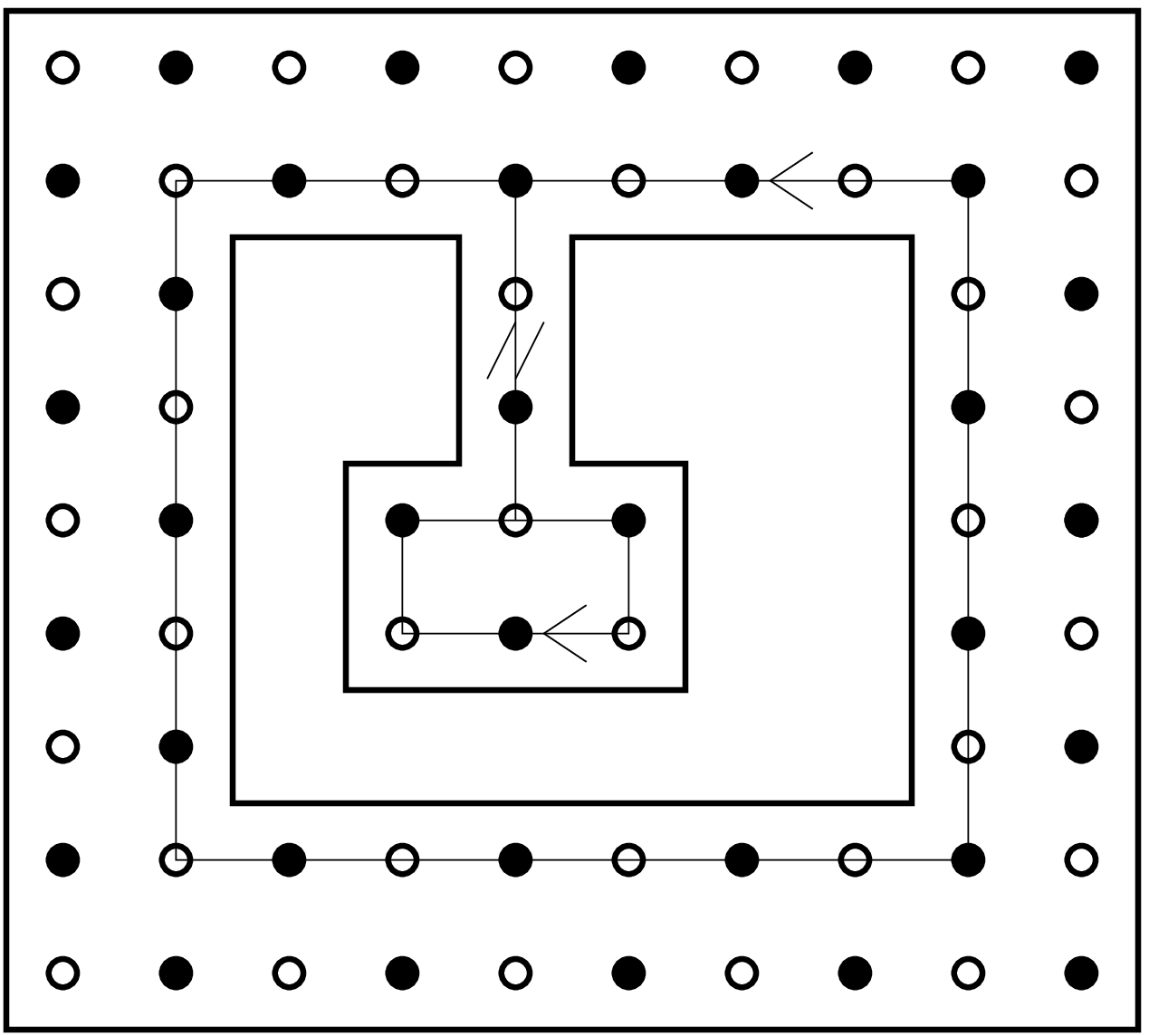} \quad
\epsfig{height=6cm,file=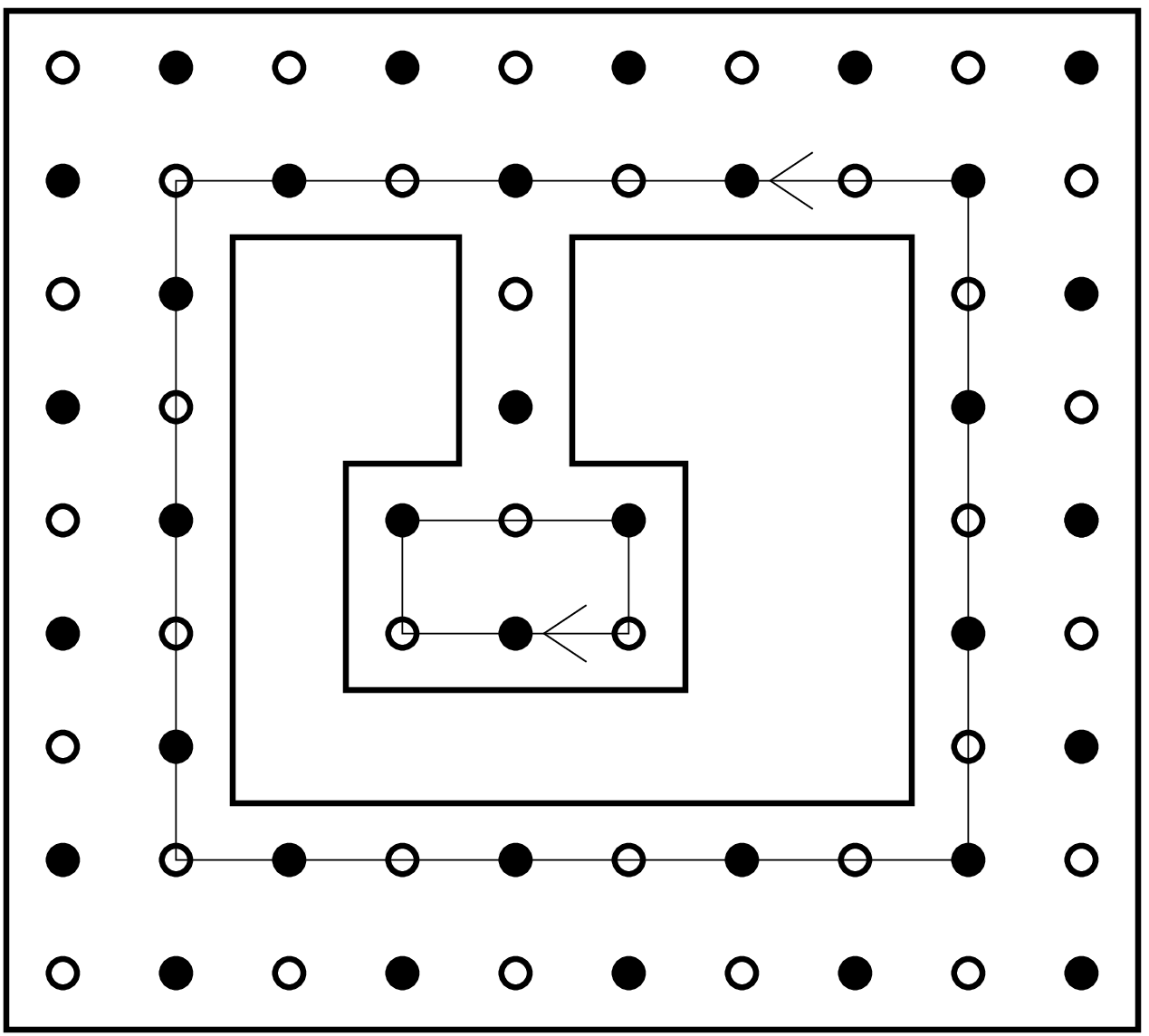}
\end{center}
\caption{The large hole}
\label{fig:badhole}
\end{figure}

There is a natural basis for $Z_1(A)$: the {\it hole basis}.
Recall that $A$ and therefore $G_A$ are embedded in $\RR^2$.
The bounded connected components of the complement of the graph $G_A$
are the {\it holes of the graph $G_A$}:
each hole gives rise to an element of the hole basis, as follows.
Many holes of the graph are bounded by four edges, the {\it small holes}.
The basis 1-cycle $s_i$ corresponding to a small hole
consists of those four edges
oriented counterclockwise (here a white-to-black edge is identified
with minus the same edge with the original black-to-white orientation).
The only remaining hole, the {\it large hole},
contains the bounded connected component of $\RR^2 - A$.
To obtain the 1-cycle $\ell$ associated to the large hole,
add the edges contained in the boundary
of the closure of the large hole, orienting them so as to have the hole
to the left.
Figure \ref{fig:badhole} below represents the same basis element $\ell$
in two ways, and is given to justify the apparently convoluted definition.

As usual, a {\it circuit} in a graph is a 2-regular subgraph;
thus, a circuit is a simple closed curve contained in $G_A$.
Oriented circuits are 1-cycles and generate $Z_1(A)$:
each $s_i$ is an oriented circuit
and $\ell$ is a sum of disjoint oriented circuits.
The decomposition in the hole basis of a circuit $c$
oriented counterclockwise is very geometric.
From the region of the plane surrounded by $c$,
remove $G_A$ to obtain a disjoint union of holes $h_1, \ldots, h_F$:
we clearly have $c = \sum h_i$.

We define three group homomorphisms by prescribing their values
on the hole basis.

\begin{defin}
\label{defin:homohole}
On the hole basis,
the {\it flux} $\phi: Z_1(A) \to \ZZ$,
{\it volume} $\nu: Z_1(A) \to \ZZ$ and
{\it sign} $\sigma: Z_1(A) \to \{ \pm 1 \}$
have the values
\begin{itemize}
\item{$\phi(s_i) = 0$, $\phi(\ell) = 1$;}
\item{$\nu(s_i) = 1$, $\nu(\ell) = 0$;}
\item{$\sigma(s_i) = -1$, $\sigma(\ell) = (-1)^{k+1}$,
where the boundary of the large hole has length $2k$.}
\end{itemize}
\end{defin}

Given the example in Figure \ref{fig:badhole},
the concept of length of $\ell$ may appear confusing.
Consider the boundary of the large hole:
this consists of $k_{\text{simple}}$ edges
which bound the large hole on one side only
and of $k_{\text{double}}$ edges
which bound the large hole on both sides
(in Figure \ref{fig:badhole}, $k_{\text{simple}} = 20$
and $k_{\text{double}} = 3$).
Let $2k = k_{\text{simple}} + 2 k_{\text{double}}$.

The definitions of flux and sign of a tiling given in the introduction
are special cases of the new definitions,
as will be seen in Corollary \ref{coro:newdefs}.
For the moment, we provide geometric meaning to these homomorphisms.

\begin{prop}
For any oriented circuit $c$,
\begin{enumerate}[\rm (a)]
\item{the flux $\phi(c)$ is the element of $H_1(A) = \ZZ$
associated with $c$, i.e.,
$\phi(c) = -1$, $1$ or $0$ depending whether $c$
surrounds the annulus's hole clockwise, counterclockwise or not at all; }
\item{the volume $\nu(c)$ is the number of small holes 
surrounded by $c$,
with a negative sign if $c$ is oriented clockwise; }
\item{the sign $\sigma(c)$ equals $(-1)^{b+w+k+1}$,
where $2k$ is the length of $c$ and
$b$ and $w$ are the number of black and white vertices
in the interior of $c$. }
\end{enumerate}
\end{prop}

{\nobf Proof: }
Let $c$ be an arbitrary oriented circuit.
For flux and volume, the proposition follows directly from
the decomposition of $c = h_1 + \cdots + h_F$ in the hole basis;
we give the details for sign.
Let $b$, $w$ and $k$ be as above.
By definition, $\sigma(c) = \prod \sigma(h_i)$ and
$$ \sigma(c) = (-1)^{F + \sum_i k_i}, $$
where $2 k_i$ is the length of the boundary of $h_i$.
Notice that the total number of edges on or inside $c$ is
$E = k + \sum_i k_i$.
Also, the number of vertices of $\cG$ on or inside $c$ is
$V = 2k + b + w$.
Since the Euler characteristic of the disk surrounded by $c$ is 1,
$F - E + V = 1$ and
$F = E - V + 1 = 1 - k - b - w + \sum_i k_i$, whence
$$\sigma(c) = (-1)^{b+w+k+1}.$$
\qed

So far, $\phi$, $\sigma$ and $\nu$ have been defined in $Z_1(A)$.
However, 
in the introduction, we defined $\phi(t; t_0)$ and $\sigma(t; t_0)$
for a tiling $t$ relative to a base tiling $t_0$.
We now relate this old definition of $\phi$ and $\sigma$ with
the new one.
Notice that a tiling $t$ does not belong to $Z_1(A)$
but the difference $t - t_0$ does.

\begin{coro}
\label{coro:newdefs}
Let $t$ be an arbitrary tiling of a quadriculated annulus $A$
with base tiling $t_0$. Then
$\phi(t; t_0) = \phi(t-t_0)$ and
$\sigma(t; t_0) = \sigma(t-t_0)$.
\end{coro}

{\nobf Proof: }
This follows from the proposition keeping in mind that the difference
between two tilings is a sum of disjoint oriented circuits $c_i$
surrounding tiled regions, so that $b_i = w_i$
and $\sigma(c_i) = (-1)^{k_i + 1}$.
\qed

In parallel with this corollary, we also set $\nu(t; t_0) = \nu(t - t_0)$,
the {\it volume} of a tiling $t$ relative to a base tiling $t_0$.
Also, the reader may check that $\phi(t; t_0) = \phi(t, \xi) - \phi(t_0, \xi)$
for an arbitrary cut $\xi$.

Let $\cA$ a commutative algebra with unit over $\CC$;
in this paper, $\cA$ is typically a ring of Laurent polynomials
in one or two variables:
$\cA_q = \CC[q,q^{-1}]$ and $\cA_{p,q} = \CC[p,p^{-1},q,q^{-1}]$.
As usual, an attribution of invertible elements of $\cA$,
i.e., elements of the multiplicative group $\cA^\ast$, to the edges of $G_A$
is a {\it $\cA^\ast$-weight in $A$}
and turns $G_A$ into a {\it weighted graph}.
A weight $\omega$ is thus a function from the set of edges to
$\cA^\ast$ and may naturally be identified with a homomorphism
$\omega: C_1(A) \to \cA^\ast$.
In particular, for a tiling $t$,
$\omega(t)$ is the product of the weights of the edges in $t$.
An oriented circuit $c$ is the sum of its edges
which are oriented black-to-white
minus the sum of its remaining edges:
$\omega(c)$ is the product of the weights of the edges in the first class
divided by the product of the weights of the remaining edges.
For a periodic track segment $\Delta$,
we are interested in certain weights in $G_A$
for the annulus $A = \clo\Delta$.

\begin{defin}
\label{defin:kastweight}
Let $\Delta$ be a periodic track segment, $A = \clo\Delta$
and $\xi$ the induced cut.
A $\cA_{p,q}^\ast$-weight $\omega$ in $A$
is a {\em Kasteleyn weight} if:
\begin{enumerate}[\rm (a)]
\item{the weight of an edge interior to $\Delta$ belongs
to $\cA_q^\ast$,}
\item{the weight of an edge crossing $\xi$ counterclockwise
belongs to $p\cA_q^\ast$,}
\item{the weight of an edge crossing $\xi$  clockwise
belongs to $p^{-1}\cA_q^\ast$,}
\item{for any oriented circuit $c$,
$\omega(c) = \sigma(c) p^{\phi(c)} q^{\nu(c)}$.}
\end{enumerate}
\end{defin}

\begin{prop}
\label{prop:buildweight}
Any periodic track segment $\Delta$ admits Kasteleyn weights.
Furthermore, given a base tiling $t_0$ of $A = \clo\Delta$,
for any tiling $t$ of $A$ we have
$\omega(t) = \omega(t_0) \sigma(t-t_0) p^{\phi(t-t_0)} q^{\nu(t-t_0)}.$
\end{prop}

{\nobf Proof: }
Pick a maximal tree inside $G_\Delta$
and assign arbitrary weights in $\cA_q^\ast$
to the edges of this tree;
in particular, assigning 1 to all edges of the tree will work.

Given an edge $e$ in $G_A$ not in the tree, there is a unique positively 
oriented circuit $c_e$ whose edges are either $e$ or edges in the tree.
The relation $\omega(c_e) = \sigma(c_e) p^{\phi(c_e)} q^{H(c_e)}$,
which $\omega$ must satisfy from condition (d) above,
now imposes the value of the weight of $e$.
This defines a weight $\omega$ on $G_A$.

If $e$ is in $\Delta$, then $c_e$ is also in $\Delta$
and therefore $\phi(c_e) = 0$:
the weight assigned to $e$ thus belongs to $\cA_q^\ast$,
verifying condition (a).
If instead $e$ crosses $\xi$ counterclockwise (resp. clockwise),
we have $\phi(c_e) = 1$ (resp. $-1$), verifying in turn condition (b)
(resp. condition (c)).

Set $W: Z_1(A) \to \cA^\ast$ to be
\[ W(c) = \omega(c) \sigma(c) p^{-\phi(c)} q^{-\nu(c)}. \]
This function is a group homomorphism and equals 1
on all positively oriented circuits $c_e$, which generate $Z_1(A)$
(a fact which is left to the reader).
Thus, $W$ is constant equal to 1, proving that condition (d)
holds for arbitrary oriented circuits.

We can now infer the value of $\omega(t)$ for a tiling $t$:
given a fixed base tiling $t_0$, write
\[
\omega(t) = \omega(t_0) \omega(t-t_0) = 
\omega(t_0) \sigma(t - t_0) p^{\phi(t - t_0)} q^{\nu(t - t_0)} =
\omega(t_0) \sigma(t; t_0) p^{\phi(t; t_0)} q^{\nu(t; t_0)}.
\]

Clearly, all Kasteleyn weights are obtained by this construction.
\qed

Consider a periodic track segment $\Delta$ with $n$ black and $n$ white squares,
$A = \clo\Delta$, $\xi$ the induced cut
and a Kasteleyn weight $\omega$ on $G_A$.
We now construct the associated $n \times n$ {\it Kasteleyn matrix}
$M_A$ with coefficients in $\cA_{p,q}$.
Begin by separately labeling black and white squares of $\Delta$
(i.e., vertices of $G_A$);
more formally, let $\eta_b$ (resp. $\eta_w$)
be a bijection from $\{1, 2, \ldots, n\}$
to the set of black (resp. white) squares.
A pair $\eta = (\eta_b, \eta_w)$ allows us to prescribe a specific
sign to a tiling $t$:
$\sigma(t; \eta)$ is the sign of the permutation
\[ \eta_w^{-1} \circ t \circ \eta_b: \{1, \ldots, n\} \to \{1, \ldots, n\}. \]
In the formula above, $t$ is interpreted as a bijection
from black squares to white squares of $A$.
Thus, a tiling $t$ admits a sign once either a base tiling $t_0$
or a labeling $\eta$ is given:
we have $\sigma(t; t_0) = \sigma(t; \eta)/\sigma(t_0; \eta)$.

Set $(M_A)_{ij}$ to be the weight of the edge
joining $\eta_b(j)$ and $\eta_w(i)$,
or zero if no such edge exists.
From conditions (a), (b) and (c) in the definition of a Kasteleyn weight,
we may write the Kasteleyn matrix
$M_A =  p^{-1} N^-_\Delta + M_\Delta + p N^+_\Delta$
uniquely for three matrices with entries in $\cA_q$:
$(N^-_\Delta,M_\Delta,N^+_\Delta)$ is the {\it Kasteleyn triple}
associated to the Kasteleyn weight.
Non-zero entries in $N^-_\Delta$ (resp. $N^+_\Delta$)
correspond to edges crossing $\xi$ clockwise (resp. counterclockwise);
non-zero entries in $M$ correspond to edges in $\Delta$.
Notice that a zig-zag cut (or more generally a wiggling cut)
yields a null $N^-_\Delta$ or $N^+_\Delta$.

Let $\Phi_A(p,q) = \det M_A(p,q)$ be
the {\it $q$-flux polynomial} of $A$ or $\Delta$.
This notation implicitly assumes two facts:
the innocuous dependence on the choice of cut and Kasteleyn weight.
More precisely, if $A = \clo\Delta = \clo\Delta'$
and Kasteleyn weights are assigned to $\Delta$ and $\Delta'$
then $\det(p^{-1} N^-_{\Delta'} + M_{\Delta'} + p N^+_{\Delta'})$
may be obtained from
$\det(p^{-1} N^-_{\Delta} + M_{\Delta} + p N^+_{\Delta})$
by multiplication by a factor of the form $\pm p^\ast q^\ast$.
This follows from the proposition below,
which also provides a combinatorial interpretation for $\Phi_A(p,q)$.

\begin{prop}
\label{prop:qflux}
Let $\Delta$ be a periodic track segment, $A = \clo\Delta$,
fix a Kasteleyn weight $\omega$ on $G_A$, a pair of bijections $\eta$,
a base tiling $t_0 \in T_A$
and let $M_A$ be the correponding Kasteleyn matrix.
Let $\Phi_A(p,q) = \det M_A(p,q)$ be the $q$-flux polynomial of $A$.
We then have, 
$$\Phi_A(p,q) =
\sigma(t_0; \eta) \omega(t_0)
\sum_{t \in T_A} p^{\phi(t; t_0)} q^{\nu(t; t_0)}.$$
Furthermore, the special values $\Phi_A(p,1)$ and $\Phi_A(p,-1)$
agree with the flux polynomials $\Phi_A^+(p)$ and $\Phi_A^-(p)$
defined in the introduction up to a multiplicative factor of the form
$\pm p^\ast$.
\end{prop}

{\nobf Proof: }
Tilings of $A$ correspond to non-zero monomials
in the expansion of $\det M_A$.
Thus, for a fixed base tiling $t_0$,
$$\det M_A = \sum_t \sigma(t; \eta) \omega(t) =
\sigma(t_0; \eta) \omega(t_0) \sum_t p^{\phi(t - t_0)} q^{\nu(t - t_0)},
$$
the desired formula for $\Phi_A(p,q)$.
Substituting 1 and $-1$ for $q$ yields the interpretations for
$\Phi_A(p,1)$ and $\Phi_A(p,-1)$ which were called
$\Phi_A^+$ and $\Phi_A^-$ in the introduction.
\qed

\begin{figure}[ht]
\begin{center}
\epsfig{height=4cm,file=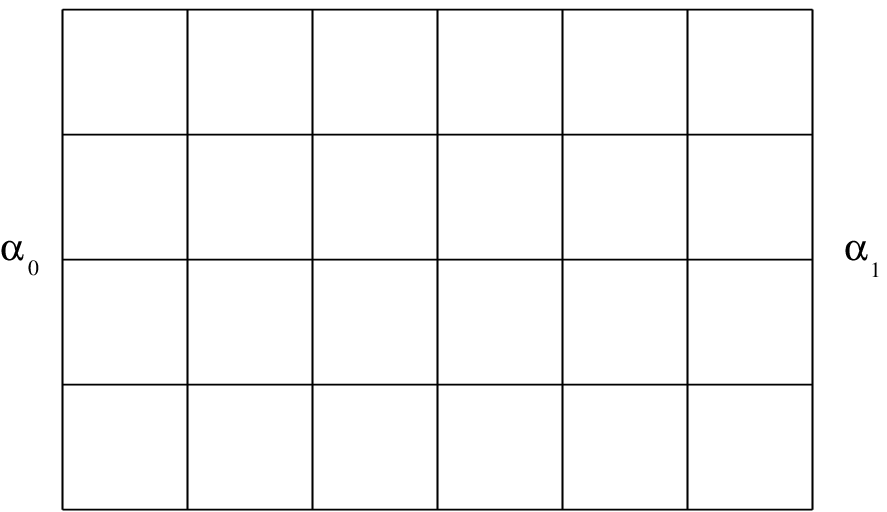}\quad
\epsfig{height=4cm,file=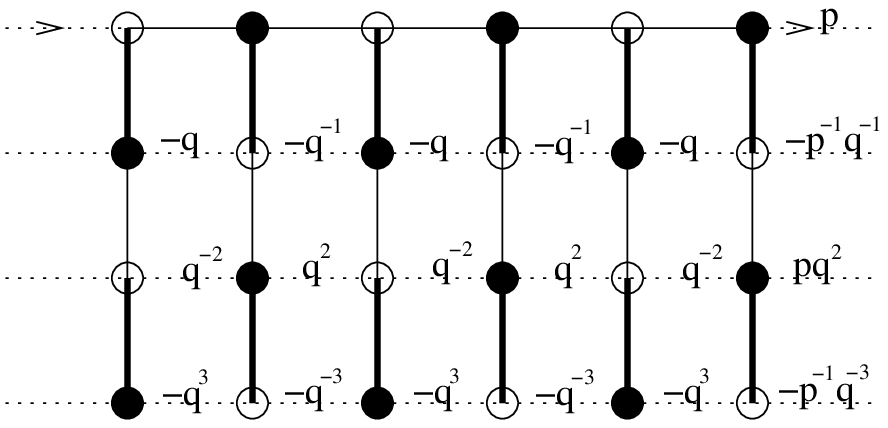}
\end{center}
\caption{Annulus and Kasteleyn weight}
\label{fig:tree}
\end{figure}

In the left of Figure \ref{fig:tree}
we show a periodic track segment $\Delta$;
the attachments are the vertical sides and the outer boundary
is the bottom side.
In the right, we show a Kasteleyn weight
constructed as in Proposition \ref{prop:buildweight}
by assigning weight 1 to the edges of the maximal tree
indicated by solid lines.
Thick lines indicate a tiling $t_0$.
Dotted lines indicate edges not belonging to the tree.
We thus obtain a $12 \times 12$ Kasteleyn matrix whose determinant is
\begin{align}
\Phi_A(p,q) &= q^{-18}\,p^{-2}\,\big(\,\,
{q}^{36}\,p^4 + \notag\\
&\phantom{}
({q}^{36}+3\,{q}^{35}+3\,{q}^{34}+4\,{q}^{33}+6\,{q}^{32}+6\,{q}^{31}+7
\,{q}^{30}+6\,{q}^{29}+6\,{q}^{28}+7\,{q}^{27}+
\notag\\
&\phantom{= +}6\,{q}^{26}+6\,{q}^{25}
+7\,{q}^{24}+6\,{q}^{23}+6\,{q}^{22}+4\,{q}^{21}+3\,{q}^{20}+3\,{q}^{
19}+{q}^{18})\,p^3 + \notag\\
&\phantom{}
({q}^{30}+3\,{q}^{29}+3\,{q}^{28}+4\,{q}^{27}+9\,{q}^{26}+12\,{q}^{25}+
16\,{q}^{24}+24\,{q}^{23}+33\,{q}^{22}+
\notag\\&\phantom{= +}
41\,{q}^{21}+45\,{q}^{20}+51\,{
q}^{19}+57\,{q}^{18}+51\,{q}^{17}+45\,{q}^{16}+41\,{q}^{15}+
\notag\\&\phantom{= +}
33\,{q}^{
14}+24\,{q}^{13}+16\,{q}^{12}+12\,{q}^{11}+9\,{q}^{10}+4\,{q}^{9}+3\,{
q}^{8}+3\,{q}^{7}+{q}^{6})\,p^2 + \notag\\
&\phantom{}
({q}^{18}+3\,{q}^{17}+3\,{q}^{16}+4\,{q}^{15}+6\,{q}^{14}+6\,{q}^{13}+7
\,{q}^{12}+6\,{q}^{11}+6\,{q}^{10}+7\,{q}^{9}+
\notag\\
&\phantom{= +}6\,{q}^{8}+6\,{q}^{7}+7
\,{q}^{6}+6\,{q}^{5}+6\,{q}^{4}+4\,{q}^{3}+3\,{q}^{2}+3\,q+1)\,p + 1\,\big)
\notag
\end{align}
Also, $\Phi_A^+(p) = p^{-2}(p^4 + 91\,p^3 + 541\,p^2 + 91\,p + 1)$
(with roots approximately equal to $-84.619$, $-6.2077$,
$-.16109$, $-.011818$) and
$\Phi_A^-(p) = p^{-2}(p^4 + p^3 + p^2 + p + 1)$,
in agreement with Theorems \ref{theo:plus} and \ref{theo:minus}.

\section{Kasteleyn matrices for covers}

In this section we present an algebraic counterpart of the juxtaposition
$\Delta^n$ of $n$ copies of a periodic track section $\Delta$,
from which we obtain some preliminary information concerning
the roots of $\Phi_A(\cdot, q)$.

\begin{prop}
Let $\Delta$ be a periodic track segment with
a Kasteleyn weight in $A = \clo(\Delta)$ with triple 
$(N^-,M,N^+) = (N^-_\Delta(q), M_\Delta(q), N^+_\Delta(q))$.
Then there is a natural Kasteleyn weight in $A^n = \clo(\Delta^n)$ with triple
$(N^-_{\Delta^n}, M_{\Delta^n}, N^+_{\Delta^n})$, where
$$N^-_{\Delta^n} = \begin{pmatrix}
\cdots & 0 &  (-1)^{n+1} N^- \\
& 0 & 0 \\
& & \vdots \end{pmatrix},
N^+_{\Delta^n} = \begin{pmatrix}
\vdots \\ 0 & 0 \\ (-1)^{n+1} N^+ & 0 & \ldots
\end{pmatrix}$$
$$M_{\Delta^n} = \begin{pmatrix}
M & N^+ & 0 & \ldots & 0 \\
N^- & M & N^+ & \ldots & 0 \\
0 & N^- & M & \ldots & 0 \\
\vdots & \vdots & \vdots & & \vdots \\
0 & 0 & 0 & \ldots & M
\end{pmatrix},$$
Thus, $M_{A^n}(p,q) =
p^{-1} N^-_{\Delta^n}(q) + M_{\Delta^n}(q) + p N^+_{\Delta^n}(q)$.
\end{prop}

{\nobf Proof:}
Here the vertices in $\Delta^n = \Delta_0 \cdots \Delta_{n-1}$
are labeled in the following order:
first those belonging to $\Delta_0$
(with the same order applied in the Kasteleyn triple for $A$),
then those in $\Delta_1$ (again with the same order) and so on.
The Kasteleyn weight in $\Delta^n$ assigns to each edge contained
in some $\Delta_k$ the same weight as in $\Delta$.
For edges trespassing a cut $\xi_{k+ \frac{1}{2}}$ clockwise
(resp. counterclockwise), $0 \le k < n-1$,
we assign as weight the corresponding entry
in $N^-_\Delta$ (resp. $N^+_\Delta$);
notice that these weights belong to $\cA_q^\ast$.
Finally, for an edge in $A^n$ trespassing
$\xi_{-\frac{1}{2}} = \xi_{n - \frac{1}{2}}$ clockwise
(resp. counterclockwise), we assign as weight $(-1)^{n+1} p^{-1}$
(resp. $(-1)^{n+1} p$) times the corresponding  entry 
in $N^-_\Delta$ (resp. $N^+_\Delta$).
This clearly yields the matrices in the statement;
we are left with proving that this construction obtains
a Kasteleyn weight in $G_{A^n}$.

Conditions (a), (b) and (c) in Definition \ref{defin:kastweight}
are clearly satisfied.
We check (d) for the hole basis.
For small holes, this is easy.
Let $\ell$ be the large hole in $G_A$ and its $n$-cover $\ell_n$
be the large hole in $G_{A^n}$.
For the Kasteleyn weight in $A$,
$\omega(\ell) = \sigma(\ell) p^{\phi(\ell)} q^{\nu(\ell)} =
(-1)^{k+1} p$, where the boundary of the large hole has length $2k$
(see Definition \ref{defin:homohole}).
We have $\sigma(\ell_n) = (-1)^{nk+1}$, $\phi(\ell_n) = 1$
and $\nu(\ell_n) = 0$ and we must check that
$\omega(\ell_n) = (-1)^{nk+1} p$, which is the case by construction.
\qed

For a Laurent polynomial $P(X) = a(X - \lambda_1)\cdots(X - \lambda_m)X^b$,
$\lambda_1,\ldots,\lambda_m \ne 0$, the $n$-th {\it root power} is
\begin{align}
P^{[n]}(X) &= a^n(X - (-1)^{n+1}\lambda_1^n)\cdots
(X - (-1)^{n+1}\lambda_m^n)X^b \notag\\
&= \prod_{j = 0, 1, \ldots, n-1} P(\zeta^{2j+1} Y), \notag
\end{align}
where $\zeta = - \exp(\pi i/n)$ and $Y^n = X$.

The reader may check that in the product above,
the terms containing $Y^k$, $k \not\equiv 0 \pmod n$, cancel out.
Let
$$\sigma_k(x_1,\ldots,x_m) =
\sum_{1 \le i_1 < \cdots < i_k \le m} x_{i_1} \cdots x_{i_k}$$
be the usual symmetric functions.
If $P(X) = (a_m X^m + \cdots + a_0)X^b =
a_m(X - \lambda_1)\cdots(X - \lambda_m)X^b$,
then clearly $P^{[n]}(X) = (a_m^{[n]} X^m + \cdots + a_0^{[n]})X^b$,
where
\begin{align}
a_k^{[n]} &= (-1)^{m-k} a_m^n 
\sigma_{m-k}((-1)^{n+1} \lambda_1^n, \ldots, (-1)^{n+1} \lambda_m^n) \notag\\
&= a_m^n 
\sigma_{m-k}((-\lambda_1)^n, \ldots, (-\lambda_m)^n). \notag
\end{align}

\begin{prop}
\label{prop:ank}
Let $\Delta$ be a periodic track segment, $A = \clo\Delta$,
$q$ a fixed complex number, $P_1(p) = \Phi_A(p,q)$
and $P_n(p) = \Phi_{A^n}(p,q)$.
Then $P_n$ is the $n$-th root power $P_1^{[n]}$.
\end{prop}

For most $q$, in particular for $q = 1$,
the number of non-zero roots of $P_1(p)$ is $f_{\max} - f_{\min}$
(recall that $f_{\max}$ (resp. $f_{\min}$) is the maximum (resp. minimum)
flow for tilings of $A$).
However, for special values of $q$ some of the extremal coefficients
of $P_1(p)$ may vanish thus causing the number of non-zero roots to go down.

{\nobf Proof:}
Let $\zeta = - \exp(\pi i/n)$ and $r^n = p$.
Let $\tilde M$ be the block diagonal matrix
with blocks
$ (\zeta^{2k+1} r)^{-1} N^-_{\Delta}(q) + M_{\Delta}(q) +
(\zeta^{2k+1} r) N^+_{\Delta}(q)$, $k = 0, 1, \ldots, n-1 $
whose determinant is clearly $P_1^{[n]}$.
The reader may check that $\tilde M = X^{-1} M_{A^n}(p,q) X$ where
$X$ is the product of the block diagonal matrix
$\diag(I, \zeta r I, \ldots, \zeta^{n-1} r^{n-1} I)$
and the discrete Fourier transform block matrix, whose $(j,k)$-block
is $\zeta^{2jk} I$, $j, k = 0, \ldots, n-1$.
\qed

The above proposition already yields considerable information
concerning the roots of the polynomials
$\Phi_A^{\pm}(p) = \Phi_A(p,\pm 1)$.
For instance, if $\lambda_1, \ldots, \lambda_m$
are the non-zero roots of $\Phi_A^+$ we must then have
$\sigma_k((-1)^{n+1} \lambda_1^n, \ldots, (-1)^{n+1} \lambda_m^n) \ge 0$
for all $k$ and $n$,
since this number counts tilings with a given flux in $\clo(\Delta^n)$
(up to a positive factor).
In order to prove the main theorems, however, we need sharper
estimates of such numbers;
these require the tools of the next sections.

\section{Tilings of track segments and $-1$-counting}

\begin{figure}[ht]
\begin{center}
\epsfig{height=3cm,file=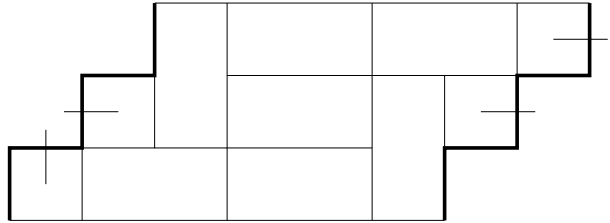}
\end{center}
\caption{Tiling of a track segment, with indices}
\label{fig:segtile}
\end{figure}

A {\it tiling} $t$ of a {\it track segment} $\Delta$ is 
a decomposition of $\Delta$ as a disjoint union of squares
with a {\it preferred edge} contained in $\az$ or $\ai$ and dominoes.
The preferred edge indicates the position of the missing half-domino.
The reader will have no difficulty in extending the tiling
of the track segment in Figure \ref{fig:segtile}
to the right but will notice that no extension to the left exists:
not every tiling of a track segment is a restriction
of a tiling of a larger region.

The {\it flux} $\phi(t, \az)$ of a tiling $t$ of a track segment $\Delta$
relative to $\az$ counts preferred edges in $\az$ with sign:
preferred edges which are sides of white (resp. black) squares
count positively (resp. negatively).
Similarly, $\phi(t, \ai)$ counts preferred edges in $\ai$
but now black and white count positively and negatively, respectively.
The difference $\phi(t, \ai) - \phi(t, \az)$ is the number of black squares
in $\Delta$ minus the number of white squares in $\Delta$
and therefore does not depend on the tiling $t$:
notice the similarity with Lemma \ref{lemma:fconst}.

We now associate a height function $\theta$
to a tiling $t$ of a track segment $\Delta$;
as before, we choose a base vertex $v_b$ in $\bo$ and set $\theta(v_b) = 0$.
Extend $\theta$ to the remaining vertices of $\Delta$ by local instructions:
it is forbidden to walk along preferred edges or
edges which are trespassed by dominoes.
The rest is similar:
$\theta$ goes up by 1 along an allowed edge
if there is a white square to its left
(or a black square to its right);
conversely, $\theta$ goes down by 1 if there is a black square to its left
(or a white square to its right).
Two height functions for the same tiling with different
base vertex differ by an additive integer constant.

For a fixed base vertex,
the restrictions of height functions to $\bo$ do not
depend on the tiling $t$.
Choose a {\it reference vertex} $v_r$ in $\bi$:
as in Lemma \ref{lemma:heightflux},
$\theta(v_r) = 4 \phi(t, \az) + c$ for some constant $c$.
The value of a height function on $\bi$ therefore
depends only on the flux of $t$.

A characterization analogous to Proposition \ref{prop:height}
is immediate.
Let $\theta$ be an integer valued function defined
on the set of vertices of $\Delta$.
Then $\theta$ is the height function of some tiling $t$
if and only if the following conditions hold:
\begin{enumerate}[\rm (a)]
\item{$\theta(v_b) = 0$;}
\item{if the oriented edge joining vertices $v_0$ to $v_1$
is in $\bi$ or $\bo$ and there is a white (resp. black)
square to its left then $\theta(v_1) - \theta(v_0) = 1$ (resp. $-1$);}
\item{if the oriented edge joining vertices $v_0$ to $v_1$
is in the interior of $\Delta$ or in an attachment and there is a white
(resp. black) square to its left then $\theta(v_1) - \theta(v_0) = 1$ or $-3$
(resp. $-1$ or $3$).}
\end{enumerate}

The {\it volume} $\nu(t)$ of a tiling $t$ of a track segment $\Delta$
is a weighted sum of values of $\theta(v)$ over the vertices of $\Delta$.
Vertices in $\bo$ and $\bi$ have weight 0 and interior vertices have weight 1.
A vertex $v$ in the interior of an attachment has weight $n/4$ if
it belongs to $n$ squares.
The volume of a track segment is thus a quarter-integer.

Tilings $t$ and $t'$ of $\Delta$ and $\Delta'$ can be juxtaposed
to produce a tiling of $\Delta\Delta'$ if and only if
$\ai$ and $\az'$ have the same shape and
the sets of preferred edges of $\ai$ and $\az'$ coincide.
The subset of preferred edges of a tiling $t$ at $\az$ (resp. $\ai$)
is the {\it 0-index} (resp. 1-index) of $t$:
thus, $t$ and $t'$ can be juxtaposed if the 1-index of $t$
coincides with the 0-index of $t'$.
More generally, a 0-index (resp. 1-index) is a subset of the set
of $n_0$ edges of $\ai$ (resp. the set of $n_1$ edges of $\az$).
Not all indices arise as the index of a tiling of $\Delta$:
for instance, the index of a tiling cannot have two preferred edges
which are sides of the same square.
Notice that a tiling of $\Delta\Delta'$
can be obtained in exactly one way as a {\it juxtaposition} $tt'$
of tilings $t$ of $\Delta$ and $t'$ of $\Delta'$.
A tiling $t$ of $\Delta$ can be {\it closed-up} to yield a tiling $\clo t$
of $\clo\Delta$ if its two indices are equal.
Either index of a tiling $t$ of $\Delta$ determines the restriction of
the height function $\theta$ to the corresponding attachment,
and therefore the flux of $t$ relative to the attachment.
Finally, a tiling $t$ of $\Delta$ and a tiling $t'$ of $\Delta'$
can be juxtaposed only if $\phi(t,\ai) = \phi(t',\az')$
(this is however not a sufficient condition).

\begin{figure}
\begin{center}
\epsfig{height=15mm,file=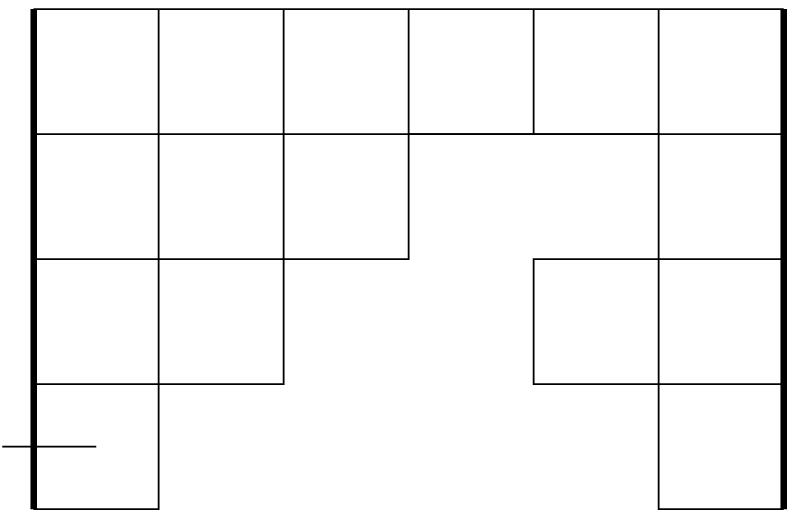}\quad
\epsfig{height=15mm,file=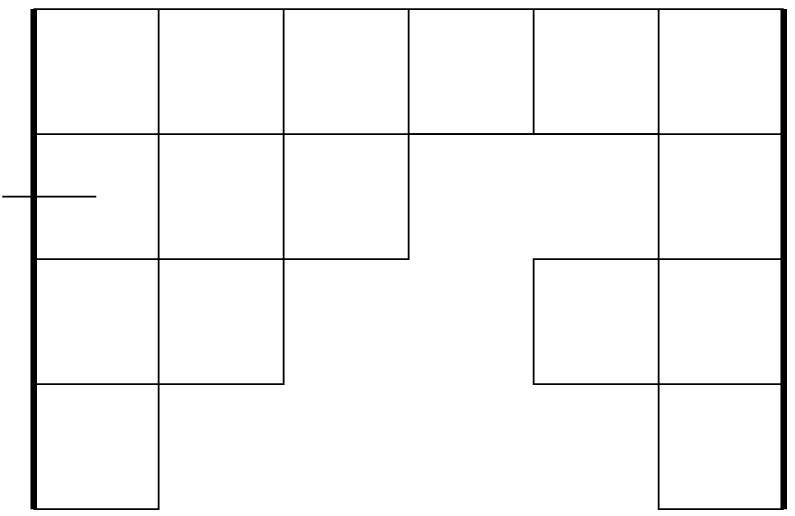}\quad
\epsfig{height=15mm,file=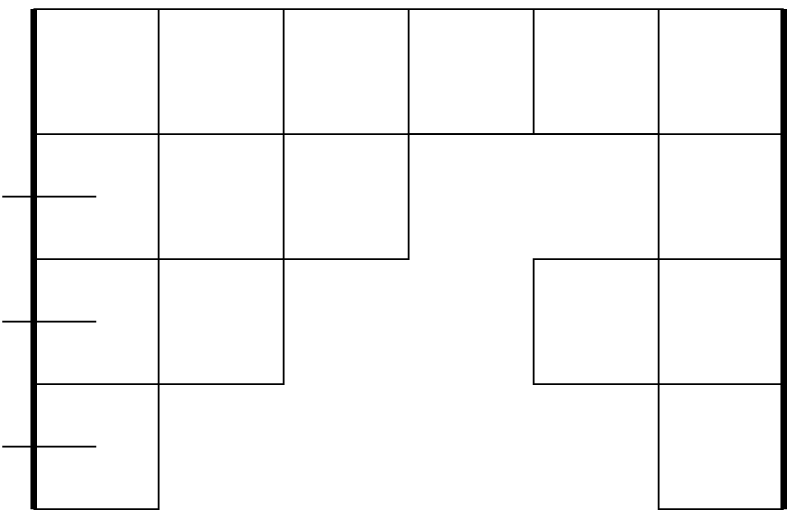}\quad
\epsfig{height=15mm,file=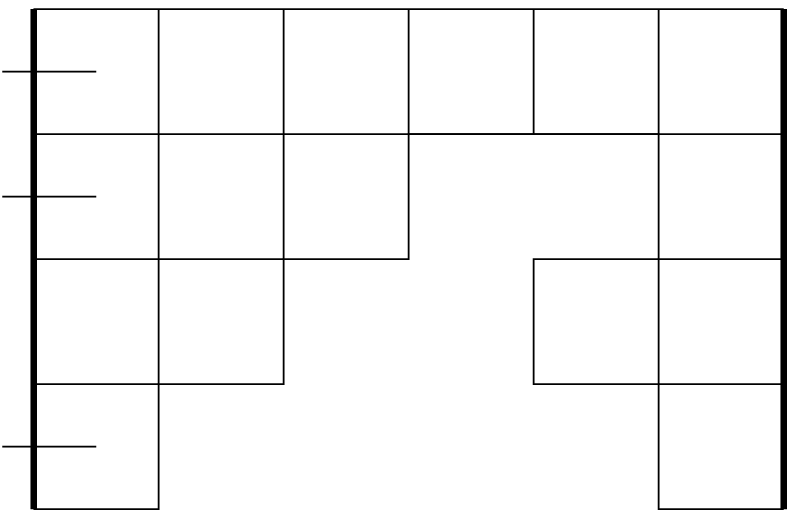}
\end{center}
\begin{center}
\epsfig{height=15mm,file=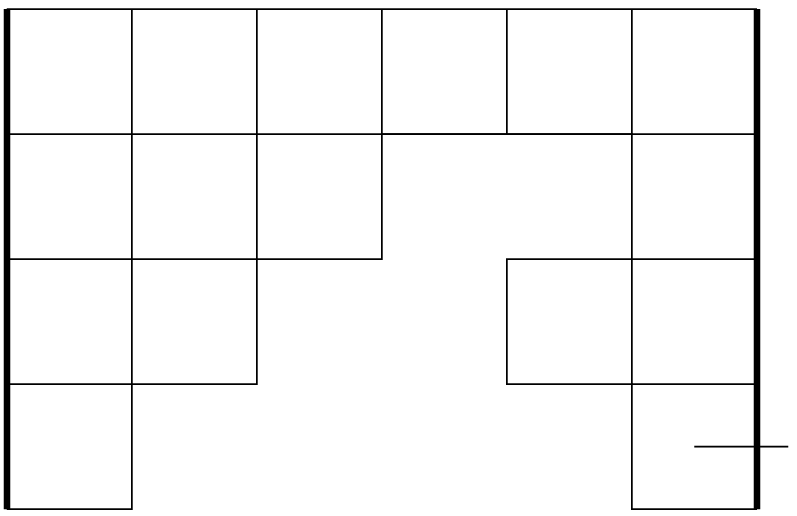}\quad
\epsfig{height=15mm,file=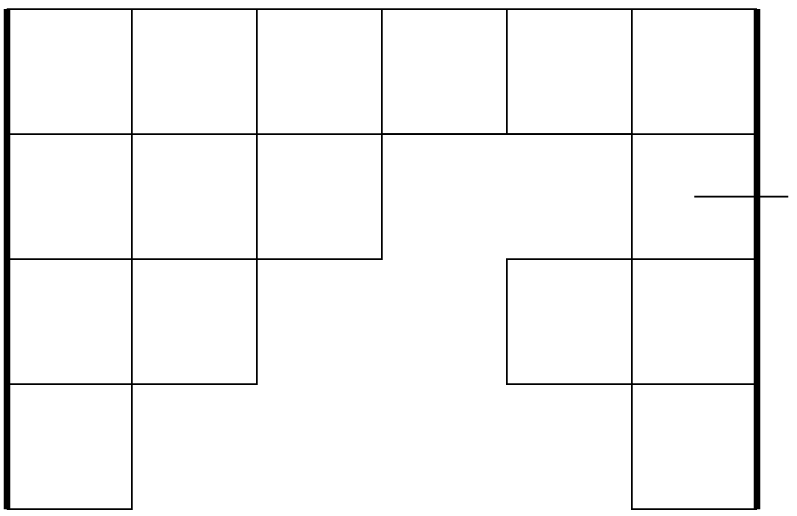}\quad
\epsfig{height=15mm,file=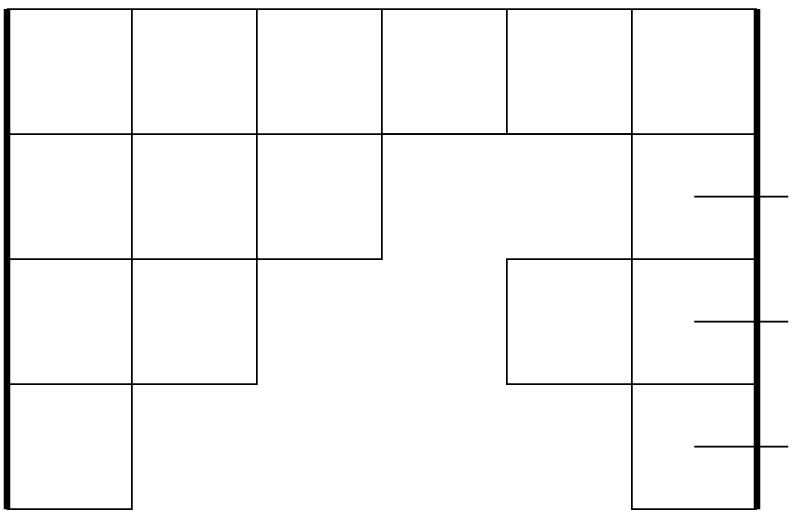}\quad
\epsfig{height=15mm,file=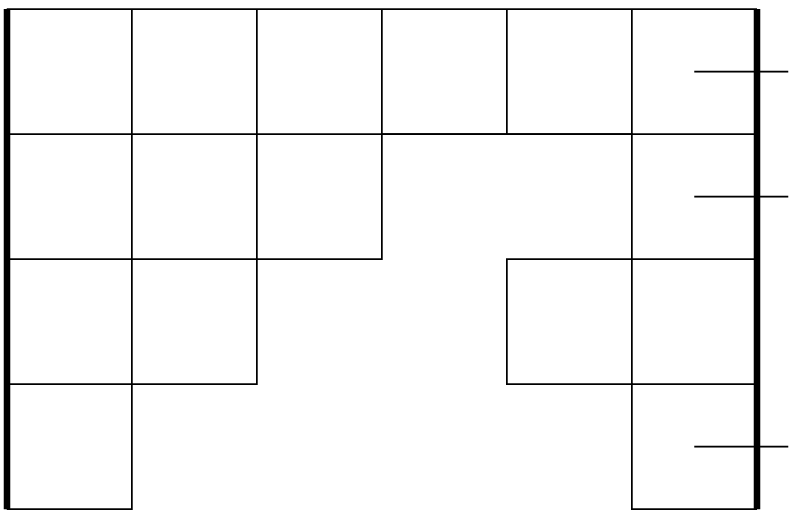}
\end{center}
\caption{Indices of flux 1}
\label{fig:nilpot}
\end{figure}

For the next definition, it is convenient to label indices.
Order the 0-indices (resp. 1-indices) so that the associated fluxes are
non-decreasing; among indices with the same flux,
use the order induced by the natural bijection
with $\{0,\ldots,2^{n_0} - 1\}$ (resp. $\{0,\ldots,2^{n_1} - 1\}$).
This assigns labels from
$\{1,\ldots,2^{n_0}\}$ and $\{1,\ldots,2^{n_1}\}$
to the sets of indices which are non-decreasing in flux.
On the first (resp. second) line of Figure \ref{fig:nilpot}
we show all 0-indices (resp. 1-indices) with associated flux $f=1$
of a track segment (the sign of the flux corresponds to painting
the lower right hand corner white).

\begin{defin}
Let $\Delta$ be a track segment with attachments $\az$ and $\ai$
of length $n_0$ and $n_1$.
The $2^{n_0} \times 2^{n_1}$ {\it connection matrix} $C_\Delta$
has its entries in $\CC[q^{1/4},q^{-1/4}]$.
The $(i_0,i_1)$ entry is the sum of $q^{\nu(t)}$
over all tilings $t$ of $\Delta$ with indices labelled $i_0$ and $i_1$.
\end{defin}

In particular, the entry $(i_0,i_1)$ is 0 if the indices labelled
$i_0$ and $i_1$ induce different values of
$\phi(\cdot,\ai)$ and $\phi(\cdot,\az')$.
The matrix $C_\Delta$ thus splits into rectangular blocks $C_{\Delta,f}$,
arranged diagonally,
labelled by the constant value $f$ of the flux on the block.

We provide a more concrete interpretation for the entries
of the connection matrix.
Given a pair of indices $(i_0,i_1)$,
we show how to {\it prune} the track segment $\Delta$
to obtain a quadriculated disk $D_{i_0,i_1}$.
Squares of $\Delta$ which are not adjacent to attachments
will be left undisturbed in $D_{i_0,i_1}$.
Consider now each square adjacent to an attachment:
if no side of this square belongs to either index,
the square also belongs to $D_{i_0,i_1}$;
if exactly one side belongs an index,
the square is not included in $D_{i_0,i_1}$;
finally, if two or more sides belong to indices,
$D_{i_0,i_1}$ is not defined.
Clearly, when $D_{i_0,i_1}$ is defined,
tilings of $D_{i_0,i_1}$ are in natural bijection
with tilings of $\Delta$ with indices $i_0$ and $i_1$;
also, when $D_{i_0,i_1}$ is undefined,
there are no tilings of $\Delta$ with indices $i_0$ and $i_1$.
The converse is not true: $D_{i_0,i_1}$ is defined even when
the corresponding fluxes are different but will not be balanced
and will therefore admit no tilings.
Summing up, the $(i_0,i_1)$ entry of $C_\Delta$ is,
up to multiplication by an integer power of $q^{1/4}$,
the $q$-counting of tilings of $D_{i_0,i_1}$
and 0 if $D_{i_0,i_1}$ is undefined.

Strictly speaking, the matrix $C_\Delta$ depends on the choice of the
base vertex. Different choices, however, merely have the effect
of multiplying the connection matrix by a power of $q^{1/4}$.

Connection matrices behave well with respect to closing-up and juxtaposition.

\begin{prop}
Let $\Delta$ be a periodic track segment
(with prescribed base and reference vertices).
Then $\Phi_A(p,q) = p^\ast q^\ast \sum_f p^f \tr C_{\Delta,f}$,
where $\Phi_A$ is the $q$-flux polynomial of the annulus $\clo\Delta$
and $p^\ast$ and $q^\ast$ denote arbitrary integer powers of
$p^{1/4}$ and $q^{1/4}$.

Let $\Delta'$ and $\Delta''$ be track segments
(with prescribed base and reference vertices) such that
the attachments $\ai'$ and $\az''$ have the same shape.
Then $C_{\Delta'\Delta''} = q^\ast C_{\Delta'} C_{\Delta''}$,
where the track segment $\Delta'\Delta''$ receives
arbitrary base and reference vertices and $q^\ast$
denotes an arbitrary integer power of $q^{1/4}$.
\end{prop}

{\nobf Proof: }
For an appropriate choice of base and reference points $v_b$ and $v_r$,
a tiling $\clo t$ of $\clo\Delta$ induces a tiling $t$ of $\Delta$
with $\nu(t; v_b) = \nu(\clo t; v_b)$,
$\phi(t; v_b, v_r) = \phi(\clo t; v_b, v_r)$
and 0- and 1-indices of $t$ are equal.
Conversely, a tiling $t$ of $\Delta$ with equal 0- and 1-indices
induces a tiling $\clo t$ of $\clo\Delta$.

For the second statement,
we may assume the base and reference vertices to belong
to the common attachment $\ai' = \az''$
since this can only modify the final result by multiplying
by an integer power of $q^{1/4}$.
Notice that with this choice of base and reference vertices,
the set of values of the flux $f$ is the same for both track segments.
The volume of a tiling $t't''$ of $\Delta'\Delta''$
clearly satisfies $\nu(t't'') = \nu(t') + \nu(t'')$.
The result now follows from the definition of matrix product.
\qed

\begin{coro}
\label{coro:trace}
Let $\Delta$ be a periodic track segment
(with prescribed base and reference vertices),
$q$ a fixed complex number and $\lambda_1, \ldots, \lambda_m$
the non-zero roots of
\begin{align}
\Phi_A(p,q) &= a p^{f_{\min}} (p^m + a_{m-1}p^{m-1} \cdots + a_0) \notag\\
&= a p^{f_{\min}} (p - \lambda_1)\cdots(p - \lambda_m). \notag
\end{align}
Then
\[ \tr C^n_{\Delta,f} =  a^n 
\sigma_{f_{\max} - f}((-\lambda_1)^n, \ldots, (-\lambda_m)^n). \]
\end{coro}

We are ready for the main two results of the section.
The matrix $C^-_{\Delta,f}$ is obtained from $C_{\Delta,f}$
by substituting $q = -1$.

\begin{theo}
\label{theo:root}
Let $\Delta$ be a periodic track segment.
The eigenvalues of the square matrices $C^-_{\Delta,f}$
are $0$ or roots of unity.
Furthermore, for any non-zero eigenvalue
the algebraic and geometric multiplicity are equal.
\end{theo}

The algebraic multiplicity of the eigenvalue $\lambda = 0$
(i.e., the multiplicity of the root zero in the characteristic
polynomial of $C^-_{\Delta,f}$)
may be different from its geometric multiplicity
(i.e., from $\dim\ker(C^-_{\Delta,f} - \lambda I)$).
For the example in Figure \ref{fig:nilpot}
(for which $f = 1$) we take as base point the central
point on the top ($\bo$) and we then have
\[
C^-_{\Delta,f} =
\begin{pmatrix}
0 & 0 & 0 & 0 \\
1 & 0 & 0 & 1 \\
1 & 0 & 0 & 1 \\
0 & 0 & 0 & 0 
\end{pmatrix}
\]
which is nilpotent. 
Nonzero entries correspond to classes of tilings
with a single element of even integer height.
Some zero entries correspond to empty classes
and other correspond to classes where cancellation occurs.

{\nobf Proof: }
Let $C = C^-_{\Delta,f}$ so that $C^n = C^-_{\Delta^n,f}$.
The entry $i_0,i_1$ of $C^n$, $n \in \{1, 2, 3, \ldots\}$,
equals, up to multiplication by an integer power of $(-1)^{1/4}$,
0 or the $-1$-counting of tilings of
the quadriculated disk $(\Delta^n)_{i_0,i_1}$.
By Theorem \ref{theo:DT} in the introduction,
these coefficients can only assume a finite number of different values.
Thus, for some distinct natural numbers $a$ and $b$ we have $C^a = C^b$,
from which the result follows.
\qed

We are now ready to prove Theorem \ref{theo:minus}.

\begin{theo}
Let $A$ be a balanced quadriculated annulus.
All non-zero roots of the signed flux polynomial $\Phi_A^-(p)$
are roots of unity. 
\end{theo}

{\nobf Proof: }
Call these roots $\lambda_1, \ldots \lambda_m$.
Let $\Delta$ be a periodic track segment with $A = \clo\Delta$.
From Theorem \ref{theo:root} we have that, for any fixed $k$, the sequence
\[ \sigma_k((-\lambda_1)^n, \ldots, (-\lambda_m)^n) \]
is eventually periodic in the variable $n$.
In other words, there exist positive integers $n_0$ and $n_1$ such that
\[ \sigma_k((-\lambda_1)^{n + Nn_1}, \ldots, (-\lambda_m)^{n + Nn_1})
= \sigma_k((-\lambda_1)^n, \ldots, (-\lambda_m)^n) \]
for all $k$, all non-negative $N$ and all $n \ge n_0$.
In particular,
\[ \sigma_k(((-\lambda_1)^{n_0n_1})^2, \ldots, ((-\lambda_m)^{n_0n_1})^2)
= \sigma_k(((-\lambda_1)^{n_0n_1}), \ldots, ((-\lambda_m)^{n_0n_1})) \]
for all $k$, meaning that taking squares effects a permutation of multiset
\[ \{ (-\lambda_1)^{n_0n_1}, \ldots, (-\lambda_m)^{n_0n_1} \}. \]
The $m!$ power of this permutation is the identity, implying
that, for all $j$,
\[ (-\lambda_j)^{2^{m!} \cdot n_0n_1} = (-\lambda_j)^{n_0n_1} \]
and thus $\lambda_j$ is a root of unity.
\qed

The reader may wonder what would have happened if instead of $\Delta$
we had chosen $\Delta'$ with $\clo\Delta = \clo\Delta' = A$.
Actually, arbitrary powers of $C_\Delta$ and $C_{\Delta'}$
have the same trace since both traces count the tilings of $A^n$.
More directly, as the reader can easily convince himself with a picture,
for sufficiently large $n$ there are
(non-periodic) track segments $\Delta_1$ and $\Delta_2$ with
$\Delta_1\Delta_2 = \Delta^n$, $\Delta_2\Delta_1 = (\Delta')^n$.

The estimate on the algebraic degree of $\lambda_j$
produced in the above proof is ridiculously large:
$\lambda_j$, being a root of an integer polynomial of degree at most
$f_{\max} - f_{\min}$ has algebraic degree at most $f_{\max} - f_{\min}$.

\section{Eventual compatibility and  $q$-counting, $q > 0$}

Two tilings $t_1$ and $t_2$ of a quadriculated surface
are said to differ by a {\it flip}
if they are identical except at a $2 \times 2$ square,
where one tiling has two horizontal
dominoes and the other has two vertical dominoes.
Two tilings $t_1$ and $t_2$ of a disk or band
differ by a flip if and only if their
height functions $\theta_1$ and $\theta_2$ differ at a single point,
the center $v$ of the square,
and then $\theta_2(v) - \theta_1(v) = \pm 4$.
Assume $\theta_2(v) - \theta_1(v) = 4$:
from the local characterization of height functions,
$v$ is a local minimum (resp. maximum) of $\theta_1$ (resp. $\theta_2$).
Moreover, flips are admissible exactly at local extrema
of a height function.

It follows from these remarks (\cite{T}, \cite{STCR}) that any two tilings
in a quadriculated disk can be joined by a finite sequence of flips.
The situation in quadriculated annuli is not so simple:
we state a special case of Theorem 4.1 in \cite{STCR}.

\begin{theo}
\label{theo:STCR}
Let $A$ be a wall-free quadriculated annulus
and let $t_1$ and $t_2$ be tilings with the same flux.
Then there exists a finite sequence of flips
joining $t_1$ and $t_2$.
\end{theo}

The hypothesis of $A$ being wall-free is necessary
as mentioned in \cite{STCR}, page 224;
the counterexample is also given in Figure \ref{fig:ladder} above.

Let $A = \clo\Delta$ be a quadriculated annulus
with universal cover $\ta = \ldots \Delta_{-1}\Delta_0\Delta_1 \ldots$.
Let $t_1$ and $t_2$ be tilings of $A$
inducing periodic tilings $\pi^{-1}(t_1)$ and $\pi^{-1}(t_2)$ of $\ta$.
A tiling $t_{1,2}$ of $\ta$ {\em interpolates} from $t_1$ to $t_2$
if there exists an integer $K$ such that $t_{1,2}$ and
$\pi^{-1}(t_1)$ (resp. $\pi^{-1}(t_2)$)
coincide in $\Delta_k$ for $k < -K$ (resp. $k > K$).
Two tilings $t_1$ and $t_2$ of $A$ are {\em eventually compatible}
if they admit interpolations $t_{1,2}$ and $t_{2,1}$.

\begin{figure}[ht]
\begin{center}
\epsfig{height=28mm,file=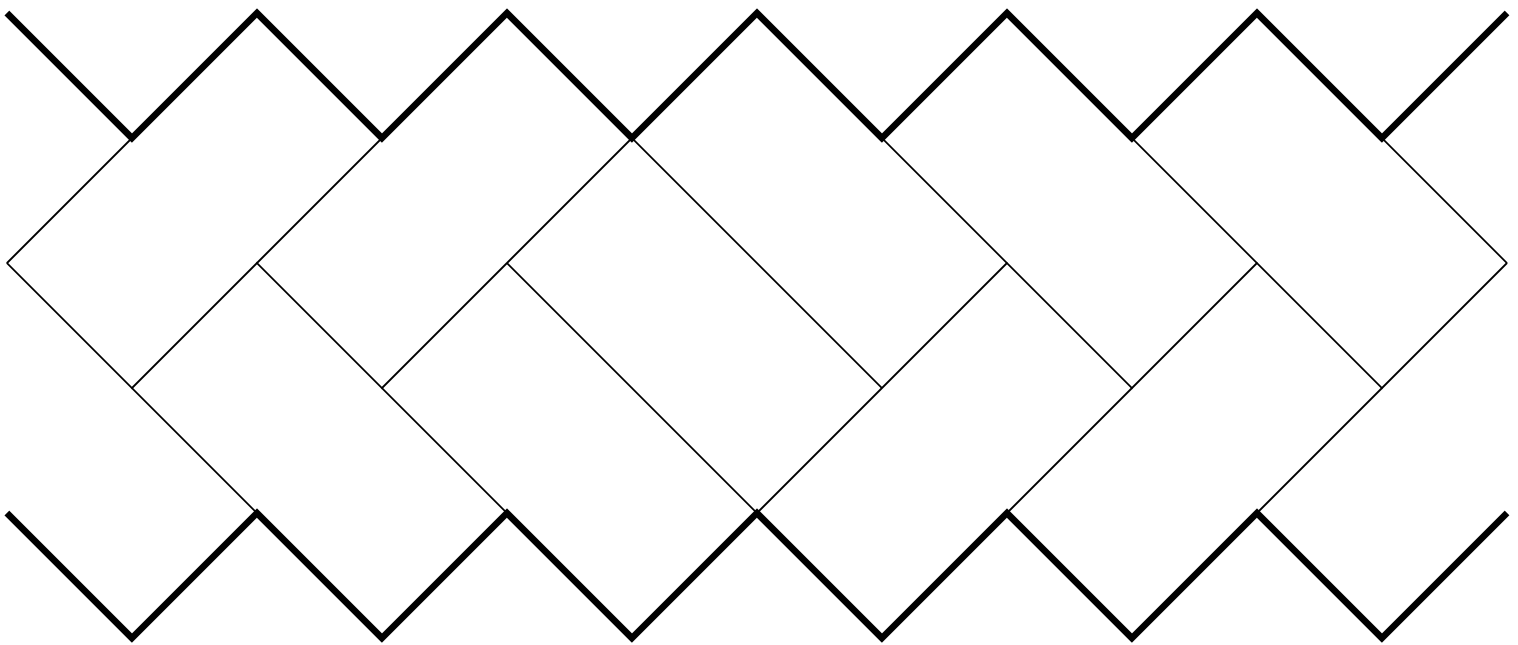} \quad
\epsfig{height=28mm,file=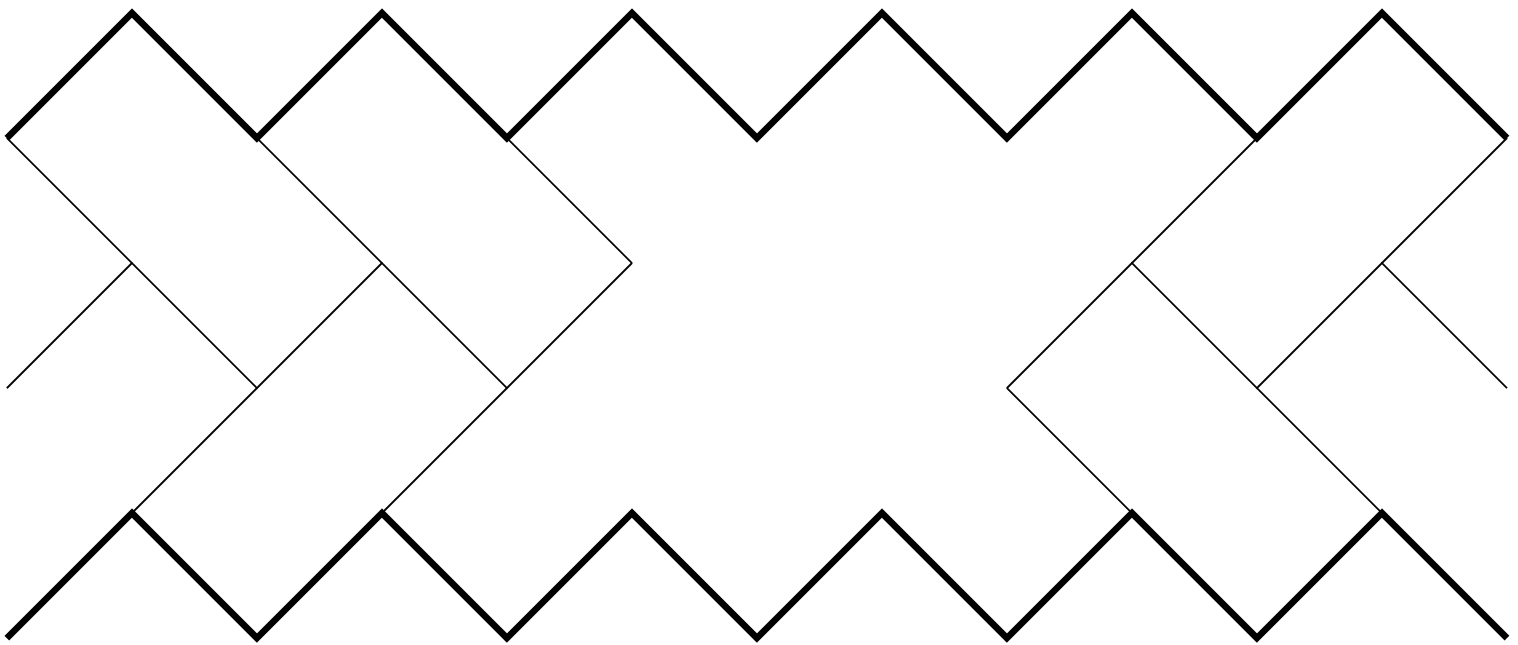}
\end{center}
\caption{Possible and impossible interpolations}
\label{fig:interpol}
\end{figure}

From Lemma \ref{lemma:fconst},
tilings with different fluxes are not eventually compatible.
Also, the tilings $t_1$ and $t_2$ in Figure \ref{fig:ladder}
have the same flux.
Figure \ref{fig:interpol} shows an interpolation $t_{1,2}$
but no interpolation $t_{2,1}$ exists (as is clear from the picture) and
hence $t_1$ and $t_2$ are not eventually compatible.
The following proposition shows that this phenomenon
does not occur in wall-free annuli.

\begin{prop}
\label{prop:compat}
Let $A$ be a balanced wall-free quadriculated annulus
and let $t_1$ and $t_2$ be tilings with the same flux
relative to the same cut.
Then $t_1$ and $t_2$ are eventually compatible.
\end{prop}

\begin{figure}[ht]
\begin{center}
\epsfig{height=3cm,file=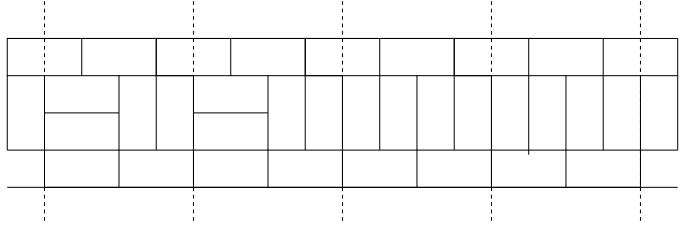}
\end{center}
\caption{Tilings differing by a flip are eventually compatible}
\label{fig:flipcompat}
\end{figure}

{\nobf Proof:}
Compatibility is an equivalence relation.
Also, tilings of $A$ differing by a flip are eventually compatible
(see Figure \ref{fig:flipcompat}).
The proposition now follows from Theorem \ref{theo:STCR}.
\qed

Let $\Delta$ be a periodic track segment and let $i$ be an index:
$i$ is {\it left-active} (resp. {\it right-active})
if there exists a tiling of the left (resp. right) half-band
$\ldots\Delta_{-1}\Delta_0$ (resp. $\Delta_1\Delta_2\ldots$)
with index $i$ at $\xi_{\frac{1}{2}}$.

An alternative definition is the following.
For a value $f$ of the flux,
we construct the {\it $f$-index graph} $\cI_f$,
an oriented graph whose vertices
are the indices of flux $f$ at $\xi_{\frac{1}{2}}$
and we join $i$ and $i'$ by an oriented edge
if there exists a tiling of $\Delta_0$
with index $i$ at $\xi_{-\frac{1}{2}}$
and $i'$ at $\xi_{\frac{1}{2}}$.
An index $i$ is left-active (resp. right-active)
if there exists an infinite path in $\cI_f$
ending (resp. starting) at $i$.

\begin{prop}
\label{prop:igraph}
Let $\Delta$ be a periodic track segment, $f$ a value of the flux
and $\cI_f$ the $f$-index graph.
Then there exists an integer $N$ such that for all $n > N$
and all indices $i$ and $i'$ the following condition holds:
there is a path in $\cI_f$
of length $n$ starting at $i$ and ending at $i'$
if and only if $i$ is  right-active and $i'$ is  left-active.
\end{prop}

{\nobf Proof:}
Let $I$ be the number of vertices in $\cI_f$
(i.e., the number of indices of flux $f$).
If there is a path of length $I+1$ starting at $i$
then the path must pass twice through some vertex
and it is straightforward to construct an eventually periodic
infinite path starting from $i$
(an infinite path $i_0i_1\ldots$ is {\it eventually periodic}
if there exist positive integers $k$ and $J$ such that
$i_{j+k} = i_j$ for all $j > J$; $k$ is the {\it period}).
Thus, if a path starts at $i$, ends at $i'$ and has length greater that $I$
then $i$ is right-active and $i'$ is left-active.

A tiling of $\clo(\Delta^n)$ induces a circuit of length $n$ in $\cI_f$;
conversely, any such circuit is induced by some tiling.
Let $i$ be a right-active and $i'$ be a left-active index.
Construct an eventually periodic path starting at $i$.
Let $k$ be its period: the path eventually correspond to a tiling
of $\clo(\Delta^k)$ which is compatible with $t_{f,\min}$
by Proposition \ref{prop:compat}
(here $t_{f,\min} \in T_\ta$ is the periodic tiling corresponding
to the minimal height function $\theta_{f,\min}$ with flux $f$ on $\ta$).
This shows us how to modify the path to obtain an eventually periodic
path of period 1, corresponding to a tiling on the right half-band
eventually coinciding with $t_{f,\min}$.
The same construction applied to the path ending at $i'$
yields a tiling of the left half-band
also eventually coinciding with $t_{f,\min}$.
It is now clear how to juxtapose these two tilings
to obtain a tiling of $\Delta_0\ldots\Delta_n$
with indices $i$ and $i'$ for any sufficiently large $n$.
\qed

An index $i$ is {\it bi-active} if both left- and right-active.
The {\it bi-active submatrix} $C_{\actv,\Delta,f}$
is obtained by removing from $C_{\Delta,f}$ 
all rows and colums associated with indices which are not bi-active.

\begin{prop}
\label{prop:eqspectra}
For any fixed complex number $q$
the spectra of $C_{\actv,\Delta,f}$ and $C_{\Delta,f}$
differ by null eigenvalues.
\end{prop}

{\nobf Proof:}
An index $i$ is
{\it strictly left-active} if left- but not right-active,
{\it strictly right-active} if right- but not left-active,
{\it inactive} if neither left- nor right-active.
We permute lines and columns of $C = C_{\Delta,f}$
so that indices appear in the following order:
left-active, bi-active, right-active, inactive,
thus turning $C$ in a $4\times 4$ block matrix
(with blocks of different sizes).
The block $C_{2,2}$ is $C_{\actv,\Delta,f}$.

We claim that $(C^n)_{2,2} = (C_{2,2})^n$
(i.e., the $2,2$ block in the $n$-th power of $C$
is the $n$-th power of its $2,2$ block).
Indeed, a monomial in the expansion of an entry
of $(C^n)_{2,2}$ is a path in $\cI_f$ starting
at a bi-active index and ending at a bi-active index:
by definition of bi-activity this path can be extended
both to the left and right and therefore,
again by definition, all indices in the path are bi-active.
In other words, the only non-zero monomials in the expansion
of $(C^n)_{2,2}$ are those in $(C_{2,2})^n$.

From Proposition \ref{prop:igraph}, for all $n > N$,
\[ C^n = \begin{pmatrix}0&0&0&0 \\ \ast & \ast &0&0 \\
\ast & \ast &0&0 \\ 0&0&0&0\end{pmatrix}. \]
It is now clear that the spectra of $C^n$ and $(C^n)_{2,2}$
differ by null eigenvalues.
From the previous paragraph, the spectra of $C^n$ and $(C_{2,2})^n$
differ by null eigenvalues and we are done.
\qed

A matrix $X$ is {\it positive} if all entries are real positive.
A square matrix $X$ is {\it eventually positive}
if $X^n$ is positive for some positive $n$.

We use the following simple modification
of the classical Perron theorem (\cite{G}).

\begin{theo}
\label{theo:Perron}
If $M$ is an eventually positive matrix
then it admits a simple positive eigenvalue $\lambda$
strictly larger than the absolute value of any other eigenvalue.
\end{theo}

\begin{theo}
\label{theo:pf}
Let $\Delta$ be a periodic track segment
with a wall-free closing-up $A = \clo\Delta$.
For any fixed positive real number $q$,
the square matrix $C_{\Delta,f}$ has a simple positive eigenvalue $\lambda$
strictly larger than the absolute value of any other eigenvalue.
\end{theo}

{\nobf Proof: }
From Proposition \ref{prop:igraph}, the bi-active submatrix
$C_{\actv,\Delta,f}$ is eventually positive.
From Theorem \ref{theo:Perron},
$C_{\actv,\Delta,f}$ admits a simple positive eigenvalue $\lambda$
strictly larger than the absolute value of any other eigenvalue.
The result now follows from Proposition \ref{prop:eqspectra}.
\qed

We are ready to restate and prove Theorem \ref{theo:qplus}.

\begin{theo}
Let $A$ be a balanced wall-free quadriculated annulus
and $q$ a fixed positive real number.
The non-zero roots (in $p$) of the flux polynomial $\Phi_A(p,q)$
are distinct real negative numbers.
\end{theo}

{\nobf Proof: }
Call $\lambda_1, \ldots \lambda_m$ the non-zero roots of $\Phi_A(p,q)$,
with $|\lambda_1| \ge |\lambda_2| \ge \ldots \ge |\lambda_m| > 0$.
Assume by induction that $\lambda_1, \ldots, \lambda_{k-1}$
are real negative and that
\[ |\lambda_1| > \cdots > |\lambda_{k-1}| > |\lambda_k| = \cdots
= |\lambda_{k'}| > |\lambda_{k'+1}|, \]
with $k' \ge k$.
We must prove that $k = k'$ and that $\lambda_k$ is real negative.

From Corollary \ref{coro:trace} and Theorem \ref{theo:pf}
we have that the sequence
of $k$-symmetric functions of powers of $\lambda_j$ satisfies
$ \sigma_k((-\lambda_1)^n, \ldots, (-\lambda_m)^n) = (1 + o(1)) B_k^n$.
In the notation of Corollary \ref{coro:trace},
the positive constant $B_k$ is the largest eigenvalue of
$C_{\Delta, f_{\max} - k}$ (for our fixed value of $q$)
divided by $a$.
The expression $ \sigma_k((-\lambda_1)^n, \ldots, (-\lambda_m)^n) $
is the sum of $k' - k + 1$ terms of the form
$((-1)^k \lambda_1\lambda_2\cdots\lambda_{k-1}\lambda_\ell)^n$,
$\ell = k, \ldots, k'$ and other terms
which grow at exponentially smaller rates.
Thus
\[ (-1)^{kn} \lambda_1^n\lambda_2^n\cdots\lambda_{k-1}^n
(\lambda_k^n + \cdots + \lambda_{k'}^n) = (1 + o(1)) B_k^n \]
whence
\[ b_k^n + \cdots + b_{k'}^n = 1 + o(1) \]
where $b_\ell = - \lambda_\ell/|\lambda_\ell|$.
Thus $b_k, \ldots, b_{k'}$ all belong to the unit circle
and we can take arbitrarily large $n$ such that 
$2 (k' - k + 1)|b_\ell^n - 1| < 1$ for all $\ell$
and therefore
\[ |(k' - k + 1) - ( b_k^n + \cdots + b_{k'}^n )| < 1/2, \]
a contradiction unless $k = k'$.
Finally, $b_k = 1$ and $\lambda_k$ is real negative.
\qed

We finish this section by estabilishing a few simple consequences
of the previous theorem.

A sequence $a_k$ of real numbers is {\it convex} (resp. {\it concave}) if
$a_{k-1} - 2a_k + a_{k+1}$ is non-negative
(resp. non-positive) for all $k$.
A sequence $a_k$ of non-negative real numbers is {\it log-concave} if
$a_k^2 \ge a_{k-1} a_{k+1}$ for all $k$.
In particular, log-concave sequences are either monotone or unimodal.

\begin{prop}
If the roots of $P(x) = a_n x^n + \cdots + a_0$
are all real negative and $a_n > 0$ then the sequence $a_k$
is log-concave.
\end{prop}

{\nobf Proof:}
We proceed by induction on the degree $n$
(the case $n = 1$ is trivial).
The result now follows from checking that if the coefficients of
$Q(x)$ form a log-concave sequence then the coefficients of
$(x + c) Q(x)$, $c > 0$ are likewise log-concave.
\qed

The converse is not true: the coefficients of $X^2 + X + 1$
form a log-concave sequence but the zeroes are not real.

\begin{coro}
Let $q$ be a fixed positive real number and
let $a_f$ be the coefficient of $p^f$ in $\Phi_A(p,q)$.
Then the sequence $a_f$ is log-concave.
\end{coro}

{\nobf Proof:}
This follows directly from Theorem \ref{theo:qplus}
and the previous proposition.
\qed

\begin{coro}
Let $b_f$ (resp. $c_f$) be the highest (resp. lowest) exponent of $q$
in non-zero terms of $\Phi_A(p,q)$ of the form $p^f q^\ast$.
Then the sequence $b_f$ (resp. $c_f$) is concave (resp. convex).
\end{coro}

{\nobf Proof:}
This follows from the previous corollary in the limit cases
$q \to 0$ and $q \to +\infty$.
\qed

{\obeylines
Nicolau C. Saldanha and Carlos Tomei
Depto. de Matem{\'a}tica, PUC-Rio
R. Mq. de S. Vicente 225
Rio de Janeiro, RJ 22453-900, Brazil
nicolau@mat.puc-rio.br
tomei@mat.puc-rio.br
}

\end{document}